# ASYMPTOTIC ERROR FOR THE MILSTEIN SCHEME FOR SDEs DRIVEN BY CONTINUOUS SEMIMARTINGALES

BY LIQING YAN

*University of Florida*

A Milstein-type scheme was proposed to improve the rate of convergence of its approximation of the solution to a stochastic differential equation driven by a vector of continuous semimartingales. A necessary and sufficient condition was provided for this rate to be $1/n$ when the SDE is driven by a vector of continuous local martingales, or continuous semimartingales under an additional assumption on their finite variation part. The asymptotic behavior (weak convergence) of the normalized error processes was also studied.

**1. Introduction.** We consider a general $q$-dimensional stochastic differential equation (SDE) driven by a vector of continuous semimartingales $Y \in \mathbb{R}^d$ on time interval $[0,1]$ with the starting point $x_0 \in \mathbb{R}^q$,

$$X_t = x_0 + \int_0^t f(X_s)\,dY_s, \tag{1}$$

where $f$ denotes a matrix of functions from $\mathbb{R}^q$ into $\mathbb{R}^q \otimes \mathbb{R}^d$. This equation includes the classical Itô-type SDE:

$$X_t = x_0 + \int_0^t a(X_s)\,dW_s + \int_0^t b(X_s)\,ds, \tag{2}$$

with $a, b$ matrices of functions and $W$ a multidimensional Brownian motion. We refer to [10] for the stochastic integral with respect to a semimartingale. In applications, one wants to find the expectation of some functional of the solution of the SDE (1), for example, such quantities can be the prices of financial derivatives, such as options. Due to the simulation difficulties, one usually combines a numerical solution $\tilde{X}^n$, with a Monte Carlo technique to approximate the expectation. For example, one can use $\mathbb{E}[h(\tilde{X}_1^n)]$ to approximate $\mathbb{E}[h(X_1)]$ and use the Riemann summation

$$\mathbb{E}\left[\frac{1}{n}\sum_{i=1}^{n} h(\tilde{X}_{t_i}^n)\right] \text{ to approximate } \mathbb{E}\left[\int_0^1 h(X_s)\,ds\right].$$

Of cause, $\{\mathbb{E}[h(\tilde{X}_{t_i}^n)], i=1,\dots,n\}$ are unknown, however, Monte Carlo method can be used to estimate them by simulating the paths of the numerical solution $\tilde{X}^n$.

A widely used numerical method for the SDE (1) is the following continuous-type Euler scheme $\tilde{X}^n$, which is given by

$$\tilde{X}_t^n = \tilde{X}_{n(t)}^n + f(\tilde{X}_{n(t)}^n)(Y_t - Y_{n(t)}), \qquad \tilde{X}_0^n = x_0,$$

where $n(t) = k/n$, if $k/n < t \le (k+1)/n$. Many authors have studied the various convergence criteria for the Euler scheme under different cases of the matrix function $f(\cdot)$ and the driving process $Y$. For reviews, see [12] or [5]. Talay and Tubaro [13] obtained the celebrated expansion of the global error of the Euler scheme. Protter and Talay [11] studied the Euler scheme for SDEs driven by Lévy processes. Kohatsu-Higa and Protter [6] studied the Euler scheme for SDEs driven by semimartingales. Yan [14] proved the convergence of Euler scheme for SDEs without continuity assumption of $f(\cdot)$, and obtained the rate of convergence without Lipschitz conditions on $f(\cdot)$.

For SDE (2), if $a \equiv 0$, then the rate of convergence of the Euler scheme is classically $1/n$; if $a$ does not vanish, it is also classical that the rate is $1/\sqrt{n}$. However, the distribution of the normalized asymptotic error for the Euler scheme was established only recently, for SDE (2) see [7], and for SDE (1) see [2]. A necessary and sufficient condition in [2] (Theorem 1.2) was given for the rate of convergence of the Euler scheme to be $1/\sqrt{n}$ when the driving semimartingale is a continuous local martingale and $f$ is a $C^1$ function, that is, $\sqrt{n}(\tilde{X}^n - X)$ converges weakly to a process $U$, which is a solution of a known linear SDE with some additional randomness.

The rate of convergence of an algorithm certainly depends on the smoothness of $f$. If $f$ is in $C^2$, we can modify the Euler scheme to improve its rate of convergence and study its normalized asymptotic error. It is well known that (see Chapter 10 of [5]) the Milstein scheme for SDE (2) with the addition of two more terms to the Euler scheme, which is given by (for the case of $d = q = 1$), $X_0^n = 0$ and

$$\begin{aligned} X_t^n = X_{n(t)}^n &+ b(X_{n(t)}^n)(t - n(t)) + \tfrac{1}{2} b(X_{n(t)}^n) b'(X_{n(t)}^n)(t - n(t))^2 \\ &+ a(X_{n(t)}^n)(W_t - W_{n(t)}) \\ &+ \tfrac{1}{2} a(X_{n(t)}^n) a'(X_{n(t)}^n)[(W_t - W_{n(t)})^2 - (t - n(t))], \end{aligned} \tag{3}$$

increases the rate of convergence from $1/n$ to $1/n^2$ when $a \equiv 0$, and increases the rate from $1/\sqrt{n}$ to $1/n$ when $a$ does not vanish. In this paper, motivated by this fact, we give a class of SDE (1) that the Milstein scheme $X^n$ converges weakly to the solution $X$ of SDE (1) at the rate of $1/n$ and determine its asymptotic error, that is, the weak limit of $n(X_t^n - X_t)$, when $f$ is a matrix of $C^2$ functions.

Our result is of mathematical interest only, since the Milstein scheme involves stochastic integrals which cannot be simulated exactly (except, in the Brownian case, when column vectors of the diffusion coefficients, seen as $C^1$ vector fields, commute in the Lie bracket sense). From an applied point of view, the discretization error needs to be studied in the weak sense (the law of the underlying process). It is now well established that errors in pathwise sense or in $L^p$ norm lead to crude sub-optimal estimates in the weak sense. For example, for the Euler scheme $X^n$ for the solution $X$, $\sqrt{n}(X^n - X)$ converges weakly to a nonzero process (see [2]) and $n(E[f(X^n)] - E[f(X)])$ converges to a nonzero constant (see [13]). However, a recent work by Kbaier [4] uses the stable weak convergence of the normalized pathwise error in order to get a useful estimate in the weak sense for the Euler scheme for Brownian SDEs. Those technical results might help to provide estimates in the weak sense for schemes (such as the Milstein scheme), which discretize SDEs driven by general semimartingales and can be simulated.

We refer to the section of preliminaries in [2] for the weak convergence in the Skorohod topology, which is denoted by "$\Longrightarrow$." See [8] and [9] for an expository account about weak convergence of stochastic integrals. For readers' convenience, we give a definition of the $(\star)$ property for a sequence of continuous semimartingales below. See [2] for the general case.

Definition. We said that a sequence of continuous semimartingales $X^n$ satisfies $(\star)$, if $[M^n, M^n]_1 + |A^n|_1$ is tight, where $X_t^n = M_t^n + A_t^n$, $M^n$ is a local martingale and $A^n$ is a finite variation process, $[M^n, M^n]_1$ is the quadratic variation of $M^n$ over $[0,1]$, $|A^n|_1$ is the total variation of $A^n$ over $[0,1]$.

In this paper we only consider the continuous process $Y$. From [2] and this paper, especially Theorem 2.2 below, we can imagine that the error distributions in the discrete case could be much more complicated. In Section 2 we develop a fundamental result on the error distribution for the Milstein scheme for the SDEs driven by a *vector of continuous semimartingales* $Y$ (not just one Brownian motion), where we employ Kurtz and Protter's idea in [7], and overcome several technical difficulties for the Milstein scheme, a more complicated and accurate algorithm in high-dimensional cases than the Euler scheme. Various convergence rates are obtained for three cases:

(i) $Y$ is a continuous process of finite variation in Section 3,
(ii) $Y$ is a continuous local martingale in Section 4 and
(iii) $Y$ is a continuous semimartingale in Section 5.

An example as the application of our theory is shown in Section 6 and seven lemmas are given in Section 7.

**2. Results on the error distribution.** Let $Y=(Y^i)_{1\le i\le d}$ be a vector of continuous semimartingales on a stochastic basis $(\Omega,\mathcal{F},\mathcal{F}_t,P)$ with $Y_0=0$. We consider the $q$-dimensional SDE (1) on time interval $[0,1]$ with the starting point $x_0=(x_0^1,\dots,x_0^q)^\tau\in\mathbb{R}^q$ and $f$ being a matrix of $C^2$ functions from $\mathbb{R}^q$ into $\mathbb{R}^q\otimes\mathbb{R}^d$ with at most linear growth [i.e., $\|f(x)\|\le c(1+\|x\|)$ for some constant $c$]. Then SDE (1) has a unique strong solution $X$, see [3].

There are many numerical methods to approximate the solution $X$, such as the Euler scheme. In this paper we consider the Milstein scheme $X^n=(X^{n1},\dots,X^{nq})^\tau$, which is defined by $X_0^n=x_0$, and

$$X_t^{ni}=X_{n(t)}^{ni}+f^i(X_{n(t)}^n)(Y_t-Y_{n(t)})\\ +\operatorname{tr}\left([f(X_{n(t)}^n)]^\tau Df^i(X_{n(t)}^n)\int_{n(t)}^t(Y_s-Y_{n(s)})\,dY_s^\tau\right),$$

where $\operatorname{tr}(A)$ is the sum of the diagonal elements of matrix $A$, $f^i$ is the $i$th row vector of $f$, and $Df^i=(f_k^{ij}:=\partial f^{ij}/\partial x_k)$ is a $q\times d$ matrix, for $i=1,\dots,q$. This scheme $X_t^{ni}$ can be written as

$$X_t^{ni}=x_0^i+\int_0^t f^i(X_{n(s)}^n)\,dY_s+\int_0^t(Y_s-Y_{n(s)})^\tau[f(X_{n(s)}^n)]^\tau Df^i(X_{n(s)}^n)\,dY_s. \tag{4}$$

See (54) for an example of this Milstein scheme for SDE (2). The error process of the Milstein scheme is denoted by $U_t^n=X_t^n-X_t$. We prove that, when $f$ is a $C^2$ function and $Y$ is a continuous local martingale, $nU^n$ converges weakly to a nonzero process $U$, which satisfies a known linear SDE (8).

First we can use the same arguments as in Theorem 3.1 in [2] to show the convergence of the Milstein scheme $X^n$ to the unique solution $X$ of SDE (1).

THEOREM 2.1. *If $f$ is a $C^2$ function with at most linear growth, then $X^n$ and $X_{n(\cdot)}^n$ go to $X$ in probability.*

Before examining rates of convergence, we introduce some notation below. For a number $0\le s\le 1$, let $s^{(n)}=s-n(s)$, where $n(s)=k/n$ if $k/n<s\le$

$(k+1)/n$. For any process $V$, we write $\Delta_t^n(V) := V_t^{(n)} := V_t - V_{n(t)}$, and

$$Z_t^n := \int_0^t \Delta_s^n(Y)(dY_s)^\tau = \int_0^t Y_s^{(n)}(dY_s)^\tau, \tag{5}$$

$$M_t^{nj} := \int_0^t \Delta_s^n(Z^n)\,dY_s^j = \int_0^t \int_{n(s)}^s Y_r^{(n)}(dY_r)^\tau\,dY_s^j, \tag{6}$$

$$N_t^{nj} := \int_0^t \Delta_s^n(Y)(\Delta_s^n(Y))^\tau\,dY_s^j = \int_0^t Y_s^{(n)}(Y_s^{(n)})^\tau\,dY_s^j. \tag{7}$$

The $(j,k)$ entries of the matrix $M^{ni}$ and $N^{ni}$ are denoted by $M^{nijk}$ and $N^{nijk}$, respectively. Let $M^n = (M^{n1}, \dots, M^{nd})$, and $N^n = (N^{n1}, \dots, N^{nd})$.

THEOREM 2.2. *Let $Hf^{ij}$ be the Hessian matrix of $f^{ij}$ and $\alpha_n$ be a deterministic sequence of positive numbers. There is equivalence between the following:*

(i) *The sequence $(\alpha_n M^n, \alpha_n N^n)$ has $(\star)$ and*

$$(Y, \alpha_n M^n, \alpha_n N^n) \Longrightarrow (Y, M, N).$$

(ii) *For any starting point $x_0$ and any $C^2$ function $f$ with at most linear growth, the sequence $\alpha_n U^n$ has $(\star)$ and $(Y, \alpha_n U^n) \Rightarrow (Y, U)$.*

*In this case, the limits $M = (M^1, \dots, M^d), N = (N^1, \dots, N^d)$ and $U = (U^1, \dots, U^q)^\tau$ can be realized on the same extension of the space on which $Y$ is defined, and they are connected by*

$$\begin{aligned} U_t^i = {} & \int_0^t (U_s)^\tau Df^i(X_s)\,dY_s \\ & - \sum_{j=1}^d \sum_{k=1}^q \mathrm{tr}\Big(\int_0^t f_k^{ij}(X_s)[Df^k(X_s)]^\tau f(X_s)\,dM_s^j\Big) \\ & - \tfrac{1}{2}\sum_{j=1}^d \mathrm{tr}\Big(\int_0^t [f(X_s)]^\tau Hf^{ij}(X_s) f(X_s)\,dN^j\Big), \end{aligned} \tag{8}$$

*and $U_0^i = 0$ for $i = 1, \dots, q$.*

PROOF. First we establish a connection between $(U^n, M^n, N^n)$ through a SDE (18). Let $R_t^n = (R_t^{n1}, \dots, R_t^{nd})^\tau$, where

$$R_t^{ni} = \int_0^t [f(X_{n(s)}^n)\Delta_s^n(Y)]^\tau Df^i(X_{n(s)}^n)\,dY_s. \tag{9}$$

Let $h^i = (Df^i)^\tau f$ by (5), we have

$$R_t^{ni} = \int_0^t \Big[h^i(X_{n(s)}^n)\Delta_s^n(Y)\Big]^\tau dY_s = \operatorname{tr}\Big(\int_0^t h^i(X_{n(s)}^n)\, dZ_s^n\Big), \tag{10}$$

$$\Delta_s^n(R^{ni}) = \operatorname{tr}(h^i(X_{n(s)}^n)\Delta_s^n(Z^n)). \tag{11}$$

By the definition of the Milstein scheme in (4),

$$X_t^n = x_0 + \int_0^t f(X_{n(s)}^n)\, dY_s + R_t^n, \tag{12}$$

and therefore,

$$\Delta_s^n(X^n) = f(X_{n(s)}^n)\Delta_s^n(Y) + \Delta_s^n(R^n). \tag{13}$$

Integrating the transpose of both sides of (13) with respect to $Df^i(X_{n(s)}^n)\, dY_s$, we have

$$\int_0^t [\Delta_s^n(X^n)]^\tau Df^i(X_{n(s)}^n)\, dY_s = \int_0^t [\Delta_s^n(R^n)]^\tau Df^i(X_{n(s)}^n)\, dY_s + R_t^{ni}, \tag{14}$$

see (9) for the definition of $R_t^{ni}$. From (12) and SDE (1), we get an expression for $U_t^{ni} = X_t^{ni} - X_t^i$,

$$U_t^{ni} = \int_0^t [f^i(X_s^n) - f^i(X_s)] - [f^i(X_s^n) - f^i(X_{n(s)}^n)]\, dY_s + R_t^{ni}. \tag{15}$$

By Taylor's expansion, there exist $\bar{\xi}_s^n$ between $X_s^n$ and $X_s$, and $\xi_s^n$ between $X_s^n$ and $X_{n(s)}^n$ such that

$$f^i(X_s^n) - f^i(X_s) = (U_s^n)^\tau Df^i(\bar{\xi}_s^n), \tag{16}$$

$$f^i(X_s^n) - f^i(X_{n(s)}^n) = [\Delta_s^n(X^n)]^\tau Df^i(X_{n(s)}^n) + \tfrac{1}{2}(H_s^{ni1}, \dots, H^{nid}), \tag{17}$$

where

$$H_s^{nij} = [\Delta_s^n(X^n)]^\tau Hf^{ij}(\xi_s^n)\Delta_s^n(X^n).$$

Combining (15), (16) and (17), we have

$$\begin{aligned} U_t^{ni} &= \int_0^t (U_s^n)^\tau Df^i(\bar{\xi}_s^n)\, dY_s + R_t^{ni} \\ &\quad - \int_0^t [\Delta_s^n(X^n)]^\tau Df^i(X_{n(s)}^n)\, dY_s - \tfrac{1}{2}\sum_{j=1}^d \int_0^t H_s^{nij}\, dY_s^j. \end{aligned}$$

From (14),

$$\begin{aligned} U_t^{ni} &= \int_0^t (U_s^n)^\tau Df^i(\bar{\xi}_s^n)\, dY_s \\ &\quad - \int_0^t [\Delta_s^n(R^n)]^\tau Df^i(X_{n(s)}^n)\, dY_s - \tfrac{1}{2}\sum_{j=1}^d \int_0^t H_s^{nij}\, dY_s^j. \end{aligned}$$

By (11) and (6),

$$\int_0^t [\Delta_s^n(R^n)]^\tau Df^i(X_{n(s)}^n)\,dY_s = \sum_{j=1}^d\sum_{k=1}^q \int_0^t f_k^{ij}(X_{n(s)}^n)\Delta_s^n(R^{nk})\,dY^j$$
$$= \sum_{j=1}^d\sum_{k=1}^q \operatorname{tr}\left(\int_0^t f_k^{ij}(X_{n(s)}^n)h^k(X_{n(s)}^n)\,dM^{nj}\right).$$

The above equation about $U^n$ can be written as follows:

$$U_t^{ni} = \int_0^t (U_s^n)^\tau Df^i(\bar{\xi}_s^n)\,dY_s - \sum_{j=1}^d\sum_{k=1}^q \operatorname{tr}\left(\int_0^t f_k^{ij}(X_{n(s)}^n)h^k(X_{n(s)}^n)\,dM_s^{nj}\right)$$
$$- \tfrac{1}{2}\sum_{j=1}^d \operatorname{tr}\left(\int_0^t Hf^{ij}(\xi_s^n)\Delta_s^n(X^n)[\Delta_s^n(X^n)]^\tau\,dY_s^j\right). \tag{18}$$

Now we prove the implication (i) $\Rightarrow$ (ii). First we assume that $f, Df, Hf$ are all bounded. Since $\Delta_s^n(Y)$, $\Delta_s^n(Z^n)$ and $\Delta_s^n(R^n)$ go to zero as $n\to\infty$ [see (11)], by the definition of $M^n$ in (6) and Theorem 2.1, we have

$$\alpha_n \int_0^t \Delta_s^n(R^{ni})\,dY_s^j = \alpha_n \operatorname{tr}\left(\int_0^t h^i(X_{n(s)}^n)\,dM_s^{nj}\right) \Longrightarrow \operatorname{tr}\left(\int_0^t h^i(X_s)\,dM_s^j\right)$$

and

$$\left(\alpha_n \int_0^t \Delta_s^n(R^n)[\Delta_s^n(R^n)]^\tau\,dY_s^j, \alpha_n\int_0^t \Delta_s^n(Y)[\Delta_s^n(R^n)]^\tau\,dY_s^j\right) \Longrightarrow 0.$$

From (13) and the definition of $N^n$ in (7), we have the following identity in the sense of weak limit ($\Longrightarrow$):

$$\lim_{n\to\infty}\operatorname{tr}\left(\alpha_n\int_0^t Hf^{ij}(\xi_s^n)\Delta_s^n(X^n)[\Delta_s^n(X^n)]^\tau\,dY_s^j\right)$$
$$= \lim_{n\to\infty}\operatorname{tr}\left(\alpha_n\int_0^t Hf^{ij}(\xi_s^n)f(X_{n(s)}^n)\Delta_s^n(Y^n)[\Delta_s^n(Y^n)]^\tau f^\tau(X_{n(s)}^n)\,dY_s^j\right)$$
$$= \lim_{n\to\infty}\operatorname{tr}\left(\alpha_n\int_0^t f^\tau(X_{n(s)}^n)Hf^{ij}(\xi_s^n)f(X_{n(s)}^n)\Delta_s^n(Y^n)[\Delta_s^n(Y^n)]^\tau\,dY_s^j\right) \tag{19}$$
$$= \lim_{n\to\infty}\operatorname{tr}\left(\alpha_n\int_0^t f^\tau(X_{n(s)}^n)Hf^{ij}(\xi_s^n)f(X_{n(s)}^n)\,dN_s^{nj}\right)$$
$$= \operatorname{tr}\left(\int_0^t f^\tau(X_s)Hf^{ij}(X_s)f(X_s)\,dN_s^j\right),$$

where (19) is due to the fact that $\operatorname{tr}(AB) = \operatorname{tr}(BA)$ for matrix $A$ and $B$. And now it follows that

$$\left(Y, \alpha_n M^n, \alpha_n\operatorname{tr}\left(\int_0^t Hf^{ij}(\xi_s^n)\Delta_s^n(X^n)[\Delta_s^n(X^n)]^\tau\,dY_s^j\right)\right)_{1\le i\le q, 1\le j\le d}$$

(20)
$$\Longrightarrow \Big(Y, M, \operatorname{tr}\Big(\int_0^t f^\tau(X_s) H f^{ij}(X_s) f(X_s)\, dN_s^j\Big)\Big)_{1\le i\le q, 1\le j\le d},$$

and the sequence in the left-hand side of (20) has $(\star)$ since $(\alpha_n M^n, \alpha_n N^n)$ has $(\star)$.

A similar argument (stopping time) as in the proof of part (a) of Theorem 3.2 in [2] shows that the weak convergence in (20) and $(\star)$ property are also true for any function $f$ in (ii). Now the implication (i) $\Rightarrow$ (ii) and the equation (8) follow from the above statement and (18) through the same arguments as in Theorem 3.2 of [2].

The implication (ii) $\Rightarrow$ (i) is proved through (18) with the starting point $x_0 = (1, \ldots, 1) \in \mathbb{R}^q$ and some special matrix function $f$, which is defined below so that $M^n$ and $N^n$ can be represented as integrals with respect to $U^n$ and $Y$.

Let's consider $M^n$ first. When $d = 1$, we choose the SDE $X_t^1 = 1 + \int_0^t X_s^1\, dY_s^1$, obviously, $X^1 = \mathcal{E}(Y^1)$, which does not vanish and $X^{n1}$ does not vanish either on $[0, 1]$ for $n$ large enough. From (18), we have $dU_t^{n1} = U_t^{n1}\, dY_t^1 - X_{n(t)}^{n1}\, dM^{n111}$, which is equivalent to

$$dM_t^{n111} = U_t^{n1}/X_{n(t)}^{n1}\, dY_t^1 - 1/X_{n(t)}^{n1}\, dU^{n1}.$$

From Theorem 2.3 in [2] and Theorem 3.1, it follows that $\alpha_n M^n$ has $(\star)$ and $(Y, \alpha_n M^n) \Rightarrow (Y, M)$.

When $d = 2$, we consider the SDE $X_t^1 = 1 + a \int_0^t X_s^1\, dY_s^1 + bY_t^2$, where $a, b$ are two constants. Equation (18) becomes

$$dU_t^{n1} = aU^{n1}\, dY_t^1 - a^3 X_{n(t)}^{n1}\, dM_t^{n111} + b\, dM^{n121}. \tag{21}$$

When $a = 1, b = 0$, let $U^{n111} = U^{n1}, X^{n111} = X^{n1}$. When $a = 1, b = 1$, let $U^{n121} = U^{n1}, X^{n121} = X^{n1}$. Equation (21) leads to two equations for $M^{n111}$, $M^{n121}$:

$$dM_t^{n111} = U_t^{n111}/X_{n(t)}^{n111}\, dY_t^1 - 1/X_{n(t)}^{n111}\, dU^{n111}, \tag{22}$$

$$dM^{n121} = -dU_t^{n121} + U^{n121}\, dY_t^1 - X_{n(t)}^{n121}\, dM_t^{n111}. \tag{23}$$

For $M^{n112}, M^{n122}$, we consider the SDE in (1) with $x_0 = (1, 1)^\tau$ and

$$f = \begin{pmatrix} a & x_1 \\ x_1 & 0 \end{pmatrix}.$$

Then

$$Df^1 = \begin{pmatrix} 0 & 1 \\ 0 & 0 \end{pmatrix} \quad \text{and} \quad Df^2 = \begin{pmatrix} 1 & 0 \\ 0 & 0 \end{pmatrix},$$

and $f_1^{21}=1$, all other $f_k^{2j}$ are zeros. Letting $i=2$ in (18) we have

$$dU_t^{n2}=U^{n1}\,dY_t^1-a\,dM_t^{n112}-X_{n(t)}^{n1}\,dM_t^{n122}. \tag{24}$$

When $a=0$, $X^1=\mathcal{E}(Y^2)$, which does not vanish and $X^{n1}$ does not vanish either on $[0,1]$ for $n$ large enough. From (24), we can get an equation for $M^{n122}$ as in (22) for $M^{n111}$. When $a=1$, from (24), we can get an equation for $M^{n112}$ as in (23) for $M^{n121}$. Now we have four equations for $M^{n1}=(M^{n1ij})_{2\times 2}$. Switching the positions of $Y^1$ and $Y^2$ gives us four similar equations for $M^{n2}$. Again, from Theorem 2.3 in [2] and Theorem 3.1, it follows that $\alpha_n M^n$ has $(\star)$ and $(Y,\alpha_n M^n)\Rightarrow(Y,M)$.

When $d\geq 3$, we choose three different fixed integers $i,j,k$ between 1 and $d$, and define

$$f^{\alpha,\beta}(x_1,x_2,\dots,x_d)=x_2(I_{1i}+I_{2k})+\lambda I_{2i}+\beta I_{2j},$$

where $I_{ij}$ is a $d\times d$ matrix with its $(i,j)$ entry being 1 and other entries being zero, $\lambda$ and $\beta$ can be either 0 or 1.

Denote $X^i(f)$ the $i$th component of the solution $X$ of SDE (1) with the starting point $x_0=(1,\dots,1)\in\mathbb{R}^q$. Let $X^{kkk}=X^2(f^{0,0})$, $X^{ikk}=X^1(f^{0,0})$, $X^{iik}=X^1(f^{1,0})$, $X^{kjk}=X^2(f^{0,1})$, $X^{ijk}=X^1(f^{0,1})$. Naturally, let $X^{nijk}$ be the Milstein scheme for $X^{ijk}$ and $U^{nijk}=X^{nijk}-X^{ijk}$ for any $1\leq i,j,k\leq d$. Since $X^{kkk}=\mathcal{E}(Y^k)$, it does not vanish and $X^{nkkk}$ does not vanish either on $[0,1]$ for $n$ large enough. Letting $i=1$ or 2 in (18), with $f$ being $f^{0,0}$, $f^{1,0}$, $f^{0,1}$, respectively, we have

$$dM_t^{nkkk}=U_t^{nkkk}/X_{n(t)}^{nkkk}\,dY_t^k-1/X_{n(t)}^{nkkk}\,dU_t^{nkkk}, \tag{25}$$

$$dM_t^{nikk}=U_t^{nkkk}/X_{n(t)}^{nkkk}\,dY_t^i-1/X_{n(t)}^{nkkk}\,dU_t^{nikk}, \tag{26}$$

$$dM_t^{niik}=U_t^{nkik}\,dY_t^i-dU_t^{niik}-X_{n(t)}^{niik}\,dM_t^{nikk}, \tag{27}$$

$$dM_t^{nkjk}=U_t^{nkjk}\,dY_t^k-dU_t^{nkjk}-X_{n(t)}^{nkjk}\,dM_t^{nkkk}, \tag{28}$$

$$dM_t^{nijk}=U_t^{nkjk}\,dY_t^j-dU_t^{nijk}-X_{n(t)}^{nkjk}\,dM_t^{nikk}. \tag{29}$$

For the error processes $U^n$, under the assumption of (ii), we have that

$$(Y,\alpha_n U^{nijk},X_{n(t)}^{nijk})_{1\leq i,j,k\leq d}\Longrightarrow(Y,U^{ijk},X^{ijk})_{1\leq i,j,k\leq d},$$

and the left-hand side of the above formula has $(\star)$. From Theorem 2.3 in [2] and Theorem 3.1, it follows that

$$\alpha_n M^n \text{ has } (\star) \text{ and } (Y,\alpha U^n,\alpha_n M^n)\Longrightarrow(Y,U,M). \tag{30}$$

Next, let us consider $N^n$. Let $F(x)=\int_x^\infty e^{-y^2/2}\,dy$, for $x,y\in\mathbb{R}^1$. When $d\geq 3$, we choose $f$ as follows, for fixed $i,j,k$ between 1 and $d$:

$$f(x_1,x_2,\dots,x_d)=F(x_2)F(x_3)I_{1i}+x_2I_{2j}+x_3I_{2j}.$$

Letting $i=1$ in (18), we have

$$
\begin{aligned}
U_t^{n1} &= \int_0^t (U_s^n)^\tau Df^1(\bar{\xi}_s^n)\,dY_s - \sum_{j=1}^d \sum_{k=1}^q \operatorname{tr}\Big(\int_0^t f_k^{1j}(X_{n(s)}^n) h^k(X_{n(s)}^n)\,dM_s^{nj}\Big) \\
&\quad - \tfrac{1}{2}\sum_{j=1}^d \operatorname{tr}\Big(\int_0^t Hf^{1j}(\xi_s^n)\Delta_s^n(X^n)[\Delta_s^n(X^n)]^\tau\,dY_s^j\Big) \qquad (31)\\
&=: \mathrm{I} - \mathrm{II} - \tfrac{1}{2}\mathrm{III}.
\end{aligned}
$$

A simple calculation shows that $Df^1 = F'(x_2)F(x_3)I_{2i} + F'(x_3)F(x_2)I_{3i}$, $Df^2 = I_{2j}$, $Df^3 = I_{3k}$ and $Df^l = 0$ for $l > 3$; $f_2^{1i} = F'(x_2)F(x_3)$, $f_3^{1i} = F'(x_3)F(x_2)$ and $f_l^{1i} = 0$ for $l \neq 2$ and $l \neq 3$, and all other Hessian matrices are zeros, except

$$
Hf^{1i}(x) = F''(x_2)F(x_3)I_{22} + F'(x_2)F'(x_3)(I_{23}+I_{32}) + F(x_2)F''(x_3)I_{33}. \qquad (32)
$$

Therefore,

$$
\begin{aligned}
\mathrm{II} &= \operatorname{tr}\Big(\int_0^t f_2^{1i}(Df^2)^\tau f\,dM_s^{ni}\Big) + \operatorname{tr}\Big(\int_0^t f_3^{1i}(Df^3)^\tau f\,dM_s^{ni}\Big),\\
\mathrm{III} &= \operatorname{tr}\Big(\int_0^t Hf^{1i}(\xi_s^n)\Delta_s^n(X^n)[\Delta_s^n(X^n)]^\tau\,dY_s^i\Big).
\end{aligned}
$$

Since $(Df^2)^\tau f = I_{j2}f = X_{n(t)}^{n2} I_{jj}$, $(Df^3)^\tau f = I_{k3}f = X_{n(t)}^{n3} I_{kk}$ and $\operatorname{tr}(I_{jj}M^{ni}) = M^{nijj}$, let $F_s^{n2} = F'(X_{n(s)}^{n2})F(X_{n(s)}^{n3})X_{n(s)}^{n2}$ and $F_s^{n3} = F'(X_{n(s)}^{n3})F(X_{n(s)}^{n2})X_{n(s)}^{n3}$, then

$$
\begin{aligned}
\mathrm{II} &= \int_0^t f_2^{1i} X_{n(s)}^{n2}\,dM^{nijj} + \int_0^t f_3^{1i} X_{n(s)}^{n3}\,dM_s^{nikk} \qquad (33)\\
&= \int_0^t F_s^{n2}\,dM_s^{nijj} + \int_0^t F_s^{n3}\,dM_s^{nikk}.
\end{aligned}
$$

On the other hand, III can be split into two terms by (13), that is, $\mathrm{III} = A_t^{ni} + B_t^{ni}$,

$$
A_t^{ni} = \int [f(X_{n(s)}^n)Y_s^{(n)}]^\tau Hf^{1i}(\xi_s^n)[f(X_{n(s)}^n)Y_s^{(n)}]\,dY_s^i \qquad (34)
$$

$$
= \operatorname{tr}\int_0^t [f(X_{n(s)}^n)]^\tau Hf^{1i}(\xi_s^n)[f(X_{n(s)}^n)]\,dN_s^i, \qquad (35)
$$

$$
B_t^{ni} = \int_0^t [\Delta_s^n(R^n) + 2f(X_{n(s)}^n)Y_s^{(n)}]^\tau Hf^{1i}(\xi_s^n)\Delta_s^n(R^n)\,dY_s^i \qquad (36)
$$

$$
= \int_0^t [\Delta_s^n(R^n) + 2f(X_{n(s)}^n)Y_s^{(n)}]^\tau Hf^{1i}(\xi_s^n)\operatorname{tr}(h^i(X_{n(s)}^n)\,dM_s^{ni}) \qquad (37)
$$

$$=:\int_0^t H_s^{ni}\,\mathrm{tr}(h^i(X_{n(s)}^n)\,dM_s^{ni}),\tag{38}$$

where (37) is due to $\int_0^t \Delta_s^n(R^{ni})\,dY_s^i = \mathrm{tr}(\int_0^t h^i(X_{n(s)}^n)\,dM_s^{ni})$, which is from (11), and (35) is due to the definition of $N^{ni}$ in (7); we write

$$H_s^{ni} = [\Delta_s^n(R^n) + 2f(X_{n(s)}^n)Y_s^{(n)}]^\tau Hf^{1i}(\xi_s^n).$$

To simplify the notation, we can write (32) as $Hf^{1i}(\xi_s^n) = aI_{22} + b(I_{23} + I_{32}) + cI_{33}$, where $a,b,c$ are corresponding function values at $\xi_s^n$. Using the fact that $I_{ij}I_{kl} = \delta_{jk}I_{il}$, where $\delta_{jk} = 1$ when $j = k$, $\delta_{jk} = 0$ when $j \neq k$; to simplify the matrix product $f^\tau Hf^{1i}f$ in (35), we get

$$\begin{aligned}
A_t^{ni} &= \mathrm{tr}\left(\int_0^t (X_{n(s)}^{n2})^2 aI_{jj} + X_{n(s)}^{n2}X_{n(s)}^{n2} b(I_{jk} + I_{kj}) + (X_{n(s)}^{n2})^2 cI_{kk}\,dN_s^{ni}\right)\\
&= \int_0^t (X_{n(s)}^{n2})^2 a\,dN_s^{nijj} + 2\int_0^t X_{n(s)}^{n2}X_{n(s)}^{n2} b\,dN^{nijk} + \int_0^t (X_{n(s)}^{n2})^2 c\,dN_s^{nikk},
\end{aligned}$$

where we use the fact that $N^{ni}$ is a sysmetric matrix. Now let $\xi_s^{n2}$ and $\xi_s^{n3}$ be the second and third component of $\xi_s^n$, respectively, and let

$$\begin{aligned}
G_s^{n1} &= (X_{n(s)}^{n2})^2 a = (X_{n(s)}^{n2})^2 F''(\xi_s^{n2})F(\xi_s^{n3});\\
G_s^{n2} &= (X_{n(s)}^{n3})^2 c = (X_{n(s)}^{n3})^2 F''(\xi_s^{n3})F(\xi_s^{n2});\\
G_s^{n3} &= 2X_{n(s)}^{n2}X_{n(s)}^{n3} b = 2X_{n(s)}^{n2}X_{n(s)}^{n3} F'(\xi_s^{n2})F'(\xi_s^{n3}),
\end{aligned}$$

then $A_t^{ni}$ can be represented as integrals with respect to $N^{ni}$ as follows:

$$A_t^{ni} = \int_0^t G_s^{n1}\,dN_s^{nijj} + \int_0^t G_s^{n2}\,dN_s^{nikk} + \int_0^t G_s^{n3}\,dN^{nijk}.\tag{39}$$

Letting $F_s^{n1} = (U_s^n)^\tau Df^1(\bar{\xi}_s^n)$, combining (31), (33), (38) and (39), we get an equation for $U^n, Y, M^n, N^n$,

$$\begin{aligned}
U_t^{n1} &= \int_0^t F_s^{n1}\,dY_s - \int_0^t G_s^{n1}\,dN_s^{nijj} - \int_0^t G_s^{n2}\,dN_s^{nikk} - \int_0^t G_s^{n3}\,dN^{nijk}\\
&\quad - \int_0^t F_s^{n2}\,dM_s^{nijj} - \int_0^t F_s^{n3}\,dM_s^{nikk} - \int_0^t H_s^{ni}\,\mathrm{tr}(h^i(X_{n(s)}^n)\,dM_s^{ni}).
\end{aligned}\tag{40}$$

When $j = k$, the above equation can be written as

$$\begin{aligned}
N_t^{nijj} = \int_0^t &\{F_s^{n1}\,dY_s - dU_s^{n1}\\
&- (F_s^{n2} + F_s^{n3})\,dM_s^{nijj} - H_s^{ni}\,\mathrm{tr}(h^i(X_{n(s)}^n)\,dM_s^{ni})\}\\
&\times \{G_s^{n1} + G_s^{n2} + G_s^{n3}\}^{-1}.
\end{aligned}\tag{41}$$

Since $X^2 = \mathcal{E}(Y^2)$ and $X^3 = \mathcal{E}(Y^3)$ do not vanish, and $X^n$ converges, $G^{ni}$ also converges and does not vanish for large $n$ for $i = 1, 2, 3$. Under (ii) in Theorem 2.2, $H^{ni}, F^{ni}$ and $G^{ni}_s$ converge to finite processes. By Theorem 2.3 in [2] and Theorem 3.1, it follows from (30) that

$$\begin{aligned}&(\alpha_n U^n, \alpha_n M^n, \alpha_n N^{nijj})_{1\le i,j\le d} \text{ has } (\star) \text{ and}\\ &\qquad (Y, \alpha_n U^n, \alpha_n M^n, \alpha_n N^{nijj})_{1\le i,j\le d} \Rightarrow (Y, U, M, N^{ijj})_{1\le i,j\le d}.\end{aligned} \tag{42}$$

From (40) we can express $N^{ijk}$ as the sum of integrals with respect to $Y$, $U^n, M^n N^{nijj}$ and $N^{nikk}$. By the same arguments, it follows from (42) that

$$\begin{aligned}&(\alpha_n U^n, \alpha_n M^n, \alpha_n N^n) \text{ has } (\star) \text{ and}\\ &\qquad (Y, \alpha_n U^n, \alpha_n M^n, \alpha_n N^n) \Rightarrow (Y, U, M, N).\end{aligned} \tag{43}$$

When $d = 1$, let $f(x_1) = F(x_1)$; when $d = 2$, take

$$f(x_1, x_2) = \begin{pmatrix} F(x_1)F(x_2) & 0 \\ & x_2 \end{pmatrix} \quad \text{and} \quad f(x_1, x_2) = \begin{pmatrix} x_1 & 0 \\ 0 & F(x_1)F(x_2) \end{pmatrix};$$

we can get (43) accordingly. Now, the implication (i) $\Rightarrow$ (ii) has been proved. □

**3. Continuous processes with finite variation.** When the driving process $Y$ in (1) is of finite variation and continuous, SDE (1) is an "$\omega$-wise" (deterministic) equation. Without loss of generosity, we assume that $Y$ is a nondecreasing function of $t$ in this section. It is classical that the rate of convergence of the Milstein scheme for an ordinary differential equations (ODE) is $1/n^2$. Here we give a necessary and sufficient condition for the rate being $1/n^2$ and its error distribution, which are useful for our main results in the next section when the driving process $Y$ in (1) is a continuous local martingale.

To simplify our notation, first we assume that $Y \in \mathbb{R}^1$ and let

$$N^n_t := N^n_t(Y) := \int_0^t (Y_s - Y_{n(s)})^2 \, dY_s. \tag{44}$$

An integration by parts shows that

$$\begin{aligned}3N^n_t(Y) &= \sum_{i=1}^{[nt]} (Y_{i/n} - Y_{(i-1)/n})^3 + (Y_t - Y_{n(t)})^3\\ &= \sum_{i=1}^{n} (Y_{(i/n)\wedge t} - Y_{((i-1)/n)\wedge t})^3.\end{aligned} \tag{45}$$

THEOREM 3.1. *Assume that $Y$ is a continuous process in $\mathbb{R}^1$ with finite variation. There is equivalence between the following:*

(a) *Equation (46) holds almost surely;*

$$Y_t = \int_0^t y_s \, ds, \qquad \int_0^1 |y_s|^3 \, ds < \infty. \tag{46}$$

(b) *The sequence of random variables $\{4^n N_1^{2^n}\}_{n\ge1}$ is tight.*
(c) *The sequence of random variables $\{n^2 N_1^n\}_{n\ge1}$ is tight.*
(d) $\sup_n n^2 N_1^n < \infty$ *a.s.*
(e) *The process $n^2 N^n$ converges a.s. uniformly in time to a process $N$.*

*In this case we have $N_t = \frac{1}{3}\int_0^t y_s^3 \, ds$ and $\sup_n n^2 \int_0^1 |dN_s^n| < \infty$ a.s.*

PROOF. The implication of (e) $\Rightarrow$ (d) $\Rightarrow$ (c) $\Rightarrow$ (b) is obvious. Let $A$ be the set of all $\omega$ such that (46) does not hold, and $B = \{\omega : 4^n N_1^{2^n} \to \infty\}$. By Lemma 3.1, we have $A \subset B$ and then (b) implies that $P(A) = 0$. So (a) follows. By Lemma 7.2, (a) implies (e) and $N_t = \frac{1}{3}\int_0^t y_s^3 \, ds$. Since $\int_0^1 |dN_s^n| \le \int_0^1 (\bar{Y}_s - \bar{Y}_{n(s)})^2 \, d\bar{Y}_s$, where $\bar{Y}_t = \int_0^t |y_s| \, ds$, and $\bar{Y}_t$ satisfies (46), by Lemma 3.1, we have $\sup_n n^2 \int_0^1 |dN_s^n| < \infty$ a.s. □

LEMMA 3.1. *Assume that $Y$ is continuous in $\mathbb{R}^1$ with finite variation and deterministic. Only two cases are possible:*

(a) *Equation (46) holds and $\sup_n n^2 N_1^n < \infty$.*
(b) *Equation (46) does not hold and $4^n N_1^{2^n} \to \infty$.*

PROOF. We can use the same notation and arguments as in the proof of Lemma 4.2 in [2]. A modification is that (4.3) in [2] is replaced by $n^2 N_1^n = \int (L_n)^3 \, d\lambda = \int (L_n)^3 V_n \, d\nu$, where $L_n(s) = n(Y_{i/n} - Y_{(i-1)/n})$ for $s \in ((i-1)/n, i/n]$ (by the way, in (4.3) of [2] $\mu'$ should be replaced by $\nu$). □

Now consider the case that $Y \in \mathbb{R}^d$. Recalling the notation of $N^{nijk}$ and $M^{nijk}$ in (6) and (7), we have the following theorem.

THEOREM 3.2. *Assume that $Y = (Y^1, \dots, Y^d)^\tau$ is a continuous process in $\mathbb{R}^d$ with finite variation. For each $i = 1, \dots, d$, there is equivalence between the following:*

(a) *We have*

$$Y_t^i = \int_0^t y_s^i \, ds, \qquad \int_0^1 |y_s^i|^3 \, ds < \infty \qquad \textit{a.s.} \tag{47}$$

(b) *The sequence of random variables $\{4^n N_1^{2^n iii}\}_{n\ge1}$ is tight.*

(c) *The sequence of random variables* $\{n^2 N_1^{niii}\}_{n\geq 1}$ *is tight.*

(d) $\sup_n n^2 N_1^{niii} < \infty$ *a.s.*

(e) *The processes* $n^2 N^{niii}$ *converges a.s. uniformly in time to a process* $N^{iii}$.

*If the above cases hold for all* $i = 1, \dots, d$, *then for any* $i, j, k = 1, \dots, d$, *then the processes* $n^2 N^{nijk}$ *and* $n^2 M^{nijk}$ *converge a.s. uniformly in time* $t$ *to processes* $N^{ijk}$ *and* $M^{ijk}$, *respectively, where*

$$N_t^{ijk} = \tfrac{1}{3}\int_0^t y_s^i y_s^j y_s^k \, ds \quad \text{and} \quad M_t^{ijk} = \tfrac{1}{6}\int_0^t y_s^i y_s^j y_s^k \, ds.$$

*And we also have* $\sup_n n^2 \int_0^1 |dN_s^{nijk}| < \infty$ *and* $\sup_n n^2 \int_0^1 |dM_s^{nijk}| < \infty$ *a.s.*

Proof. By Theorem 3.1, we have that (e) $\Rightarrow$ (d) $\Rightarrow$ (c) $\Rightarrow$ (b) $\Rightarrow$ (a). By Lemma 7.2, we have (a) $\Rightarrow$ (e), and the uniform convergence in the last claim. From (59) and Lemma 7.1, it follows that

$$\begin{aligned} n^2 \int_0^1 |dM_s^{nijk}| &\leq n^2 \int_0^1 \left| \int_{n(s)}^s Y_r^{(n)j} \, dY_r^k \right| |y_s^i| \, ds \\ &\leq n^2 \int_0^1 \left( \int_{n(s)}^s |y_r^j| \, dr \int_{n(s)}^s |y_r^k| \, dr \right) |y_s^i| \, ds < \infty, \\ n^2 \int_0^1 |dN_s^{nijk}| &\leq n^2 \int_0^1 \int_{n(s)}^s |y_r^j| \, dr \int_{n(s)}^s |y_r^k| \, dr |y_s^i| \, ds < \infty. \end{aligned}$$

So we have the uniform boundedness in the last claim. $\square$

Recalling the notation that $U_t^{n1} = X_t^{n1} - X_t^1$, $X_t^{n1}$ is the first component of the Milstein scheme $X^n$ defined in (4), and $X_t^1$ is the first component of the solution $X$ of SDE (1): $X_t = x_0 + \int_0^t f(X_s) \, dY_s$. Let

$$G^{ij}(x) = Hf^{ij}(X_s)[f(X_s)]^\tau + \sum_{k=1}^d f_k^{ij}(x)[Df^k(x)]^\tau.$$

Theorem 3.3. *When* $Y = (Y^1, \dots, Y^d)^\tau$ *is a continuous process with finite variation, there is equivalence between the following:*

(a) *For any starting point* $x_0 \in \mathbb{R}^q$ *and any* $C^2$ *function* $f$ *with linear growth, the sequence of random variables* $n^2 \int_0^1 |dU_s^n|$ *is bounded a.s.*

(b) *For any starting point* $x_0$ *and all* $C^2$ *functions* $f$ *with at most linear growth, the process* $n^2 U^n$ *converges uniformly to a limit* $U$ *a.s.*

(c) *For each $i$, the process $Y^i$ has the form*

$$Y_t^i = \int_0^t y_s^i\,ds \qquad \text{with } \int_0^1 |y_s^i|^3\,ds < \infty.$$

*In this case, let $y_s = (y_s^1, \dots, y_s^d)^\tau$, then the limiting process $U = (U^1, \dots, U^q)^\tau$ is the solution of the following linear equation: $U_0 = 0$, and*

$$U_t^i = \int_0^t (U_s)^\tau Df^i(X_s) y_s\,ds - \frac{1}{6}\sum_{j=1}^d \int_0^t y_s^\tau G^{ij}(X_s) f(X_s) y_s\,ds.$$

PROOF. That (b) ⇒ (a) is obvious. By Theorems 2.2 and 3.2, we have that (c) ⇒ (b) and the last claim. The implication of (a) ⇒ (c) can be proved by the notation and results in the proof of (ii) ⇒ (i) in Theorem 2.2. First, since $X^{nkkk}$ uniformly converges to $X^{kkk} = e^{Y^k}$ (by Theorem 2.1) and $Y_t^k$ is continuous on $[0,1]$, for fixed $\omega \in \Omega$, there exist two constants $a$ and $b$ such that $0 < a < X_t^{nkkk} < b$ for $n$ large enough and all $0 \le t \le 1$, $k = 1, \dots, d$. Now from (25)–(29), $n^2 \int_0^1 |dM_s^{nijk}|$ is bounded a.s. Similarly, from (41), $n^2 \int_0^1 |dN_s^{nijj}|$ is bounded a.s. (c) follows from Theorem 3.2. □

**4. Continuous local martingales.** To simplify our notation in the proof of our main result, we first assume that the driving process $Y$ in SDE (1) is a continuous local martingale in $\mathbb{R}^1$, null at zero. With the notation of $N_t^n(Y)$ in (44), we have the following:

THEOREM 4.1. *Assume that $Y$ is a continuous local martingale in $\mathbb{R}^1$ null at zero. Let $C_t = [Y, Y]_t$ be the quadratic variation process of $Y$. There is equivalence between the following:*

(a) *There is a predictable process $c$ such that*

$$C_t = \int_0^t c_s\,ds, \qquad \int_0^1 |c_s|^3\,ds < \infty \qquad \text{a.s. for } 0 \le t \le 1. \tag{48}$$

(b) *The sequence of random variables $\{nN_1^n(Y), n \ge 1\}$ is tight.*
(c) *The sequence of random variables $\{n^2[N^n(Y), N^n(Y)]_1; n \ge 1\}$ is tight.*
(d) *The sequence of random variables $\{n^2 N_1^n(C), n \ge 1\}$ is tight.*
(e) $\sup_{n \ge 1} n^2 N_1^n(C) < \infty$, *a.s.*

PROOF. The equivalence between (a), (d) and (e) can be deduced from Theorem 3.1. The equivalence between (b) and (c) comes from the fact that, for a sequence of continuous local martingales $M_t^n$, $\{M_1^n, n \ge 1\}$ is tight if and only if $\{[M^n, M^n]_1, n \ge 1\}$ is tight. This is because there are two universal constants $a_1$ and $a_2$ such that $a_1 \mathbb{E}([M^n, M^n]^T) \le (\mathbb{E} M_T^n)^2 \le a_2 \mathbb{E}([M^n, M^n]^T)$ for any stopping time $T$. (a) ⇒ (b) is obvious after we prove the weak convergence of $nN_t^{niii}$ in Theorem 4.2 under the assumption of (a), see the last paragraph of the proof of Theorem 4.2.

Next we show that (c) $\Rightarrow$ (d). Let $T(n,p)=\inf\{t:n^2[N^n(Y),N^n(Y)]_t\geq p\}$, and $T(p)=\inf\{t:C_t\geq p\}$, $T=T(p)\wedge T(n,p)$, and $\tau_i=T\wedge\frac{i}{n}$. Note that for a continuous local martingales $M_t$ null at zero, by the Itô formula, $M_t^6=6\int_0^t M_s^5\,dM_s+15\int_0^t M_s^4\,d[M,M]_s$. If $[M,M]_1$ is bounded, then $\int_0^t M_s^5 dM_s$ is a martingale, hence, $\mathbb{E}M_t^6=15\mathbb{E}\int_0^t M_s^4\,d[M,M]_s$. By the Burkholder–Gundy inequality, there is a constant $\kappa$ such that $\mathbb{E}[M,M]_t^3\leq\kappa\mathbb{E}M_t^6$. Therefore,

$$\mathbb{E}[M,M]_t^3\leq 15\kappa\mathbb{E}\int_0^t M_s^4\,d[M,M]_s.$$

Now let $M_t=Y-Y_{t\wedge\tau_{i-1}}$. Since $[M,M]_t=C-C_{t\wedge\tau_{i-1}}$,

$$\mathbb{E}(C-C_{t\wedge\tau_{i-1}})^3\leq 15\kappa\mathbb{E}\int_{t\wedge\tau_{i-1}}(Y_s-Y_{(i-1)/n})^4\,dC_s.$$

Recalling the notation of (45),

$$\begin{aligned}\mathbb{E}N_{t\wedge T}^n(C)&=\mathbb{E}\tfrac{1}{3}\sum_{i=1}^n(C-C_{t\wedge\tau_{i-1}})^3\\&\leq 5\kappa\mathbb{E}\int_0^{t\wedge T}(Y_s-Y_{(i-1)/n})^4\,dC_s=5\kappa\mathbb{E}[N^n(Y),N^n(Y)]_{t\wedge T}.\end{aligned}$$

Let $t=1$ in the above inequality, then $\mathbb{E}[n^2N_{1\wedge T}^n(C)]\leq 5\kappa p$. Now for any $q>0$,

$$P(n^2N_1^n(C)>q)\leq P(n^2N_{1\wedge T}^n(C)>q)+P(T<1)\leq 5\kappa p/q+P(T<1).$$

Since $P(T<1)$ goes to 0 uniformly in $n$ as $p\to\infty$ by (c), we have

$$\lim_{p\to\infty}\sup_n P(n^2N_1^n(C)>q)=0,$$

which implies (d). $\square$

Now we assume that the driving process $Y=(Y^1,\ldots,Y^d)^\tau$ in SDE (1) is a continuous $d$-dimensional local martingale on $(\Omega,\mathcal{F},\mathcal{F}_t,P)$, null at zero. Let $C_t=(C_t^{ij})_{1\leq i,j\leq d}$ be the quadratic variation processes of $Y$, that is, $C_t^{ij}=[Y^i,Y^j]_t$. Our main result is as follows.

THEOREM 4.2. *Assume that $Y$ is a continuous local martingale in $\mathbb{R}^d$ null at zero. There is equivalence between the following:*

(a) *There is a $d\times d$ nonnegative matrix-valued predictable process $c$ such that*

$$C_t=\int_0^t c_s\,ds,\qquad\int_0^1\|c_s\|^3\,ds<\infty\qquad\textit{for } 0\leq t\leq 1. \tag{49}$$

(b) *For each $i$, the sequence of random variables (r.v.) $\{nN_1^{niii}, n \geq 1\}$ is tight.*

(c) *For each $i$, the sequence of r.v. $\{n^2[N^{niii}, N^{niii}]_1, n \geq 1\}$ is tight.*
*In this case, and if $\sigma_t$ is a $d \times m$ matrix-valued process such that $c_t = \sigma_t \sigma_t^\tau$ ($\sigma^\tau$ is the transpose of $\sigma$), then there is a Brownian motion $W = (W^1, \ldots, W^m)^\tau$ on $(\tilde{\Omega}, \tilde{\mathcal{F}}, \tilde{\mathcal{F}}_t, \tilde{P})$, a possible enlarge space of $(\Omega, \mathcal{F}, \mathcal{F}_t, P)$, such that $Y_t = \int_0^t \sigma_s \, dW_s$. The sequence $(nM^{nj}, nN^{nj})_{1 \leq j \leq d}$ stably converges in law to processes $(M^j, N^j)_{1 \leq j \leq d}$, given by*

$$M_t^j = \frac{\sqrt{6}}{6} \sum_{p=1}^m \int_0^t \sigma_s \sigma_s^{jp} \, dB_s^p \sigma_s^\tau, \tag{50}$$

$$N_t^j = \frac{\sqrt{3}}{3} \sum_{p=1}^m \int_0^t \sigma_s \sigma_s^{jp} \, dV_s^p \sigma_s^\tau, \tag{51}$$

*where $B_t^p = (B_t^{pij})_{1 \leq i,j \leq m}$, $V_t^p = (V_t^{pij})_{1 \leq i,j \leq m}$ and for $1 \leq p, i, j \leq m$,*

$$V^{pij} = \sqrt{2}/2 B^{pij} + \sqrt{2}/2 B^{pji} + (\sqrt{3}/2 W^p + 1/2 \bar{W}^p) \mathbf{1}(i = j),$$

*where $\mathbf{1}(i = j) = 1$ if $i = j$, $0$ otherwise, and $\{B^{pij}, \bar{W}^p | 1 \leq p, i, j \leq m\}$ is a sequence of independent standard Brownian motions on an extension of space $(\tilde{\Omega}, \tilde{\mathcal{F}}, \tilde{\mathcal{F}}_t, \tilde{P})$ on which $Y$ and $W$ are defined and independent of $Y$ and $W$. Moreover, we have that $(nM^n, nN^n)$ has $(\star)$ and $(Y, nM^n, nN^n) \Longrightarrow (Y, M, N)$, that is, the statement* (i) *in Theorem 2.2.*

Remark 4.1. The limits of the quadratic variation of $(nM^n, nN^n, W)$ are listed in Lemma 7.6.

Proof of Theorem 4.2. By Theorem 4.1, we have the equivalence between (b) and (c). Next we show that (b) $\Rightarrow$ (a). Since $C_t - C_s$ is a non-negative matrix for any $t \geq s \geq 0$, $[C_t^{ij} - C_s^{ij}]^2 \leq [C_t^{ii} - C_s^{ii}][C_t^{jj} - C_s^{jj}]$, the Cauchy inequality implies that

$$[N_t^n(C^{ij})]^2 \leq N_t^n(C^{ii}) N_t^n(C^{jj}).$$

By the notation of $N_t^{niii}$ in (7) and $N_t^n(Y^i)$ in (44), we have $N_t^{niii} = N_t^n(Y^i)$. From the implication of (b) $\Rightarrow$ (e) in Theorem 4.1, $\sup_n n^2 N_1^n(C^{ii}) < \infty$ for all $i$. Therefore, $\sup_n n^2 N_1^n(C^{ij}) < \infty$ for all $i$ and $j$. By Theorem 3.1, we conclude (a).

It remains to prove that (a) implies the last claims, since (b) or (c) follows from the last claims. By Lemma 7.6 and Theorem 4.1 of [1], the sequence $(nM^{np}, nN^{np})_{1 \leq p \leq d}$ *stably* converges in law to the process $(M^p, N^p)_{1 \leq p \leq d}$ of (50) and (51). By the same arguments as [2], it follows that the triple $(Y, nM^{np}, nN^{np})_{1 \leq p \leq d}$ converges in law to $(Y, M^p, N^p)_{1 \leq p \leq d}$ for the product topology on $\mathbb{D}(\mathbb{R}^d) \times \mathbb{D}(\mathbb{R}^{2d^3})$, and since all these processes are continuous, we also have convergence for the Skorohod topology on $\mathbb{D}(\mathbb{R}^{d+2d^3})$. Finally, Lemma 7.6 implies $(\star)$ for the sequence $(nM^n, nN^n)$. □

THEOREM 4.3. *For a continuous local martingale* $Y = (Y^1, \dots, Y^d)^\tau$ *with quadratic variation* $C_t = [Y, Y]_t$, *there is equivalence between the following:*

(a) *For any starting point* $x_0 \in \mathbb{R}^q$ *and any* $C^2$ *function* $f$ *with linear growth, the sequence of random variables* $nU^{n*}$ *is tight, where* $X^*$ *means that* $\sup_{0 \le t \le 1} |X_t|$.

(b) *For any starting point* $x_0$ *and any* $C^2$ *functions* $f$ *with at most linear growth, the sequence* $(Y, nU^n)$ *converges weakly to limit* $(Y, U)$.

(c) *We have* (49).

*In this case, the process* $U = (U^1, \dots, U^q)^\tau$ *is the solution of the linear equation* (8), *with* $M, N$ *given in* (50) *and* (51).

PROOF. That (b) ⇒ (a) is obvious. From Theorems 2.2 and 4.2, we have that (c) ⇒ (b) and the last claim. The implication of (a) ⇒ (c) can be proved by the notation and results in the proof of (ii) ⇒ (i) in Theorem 2.2. First, from (25)–(29) and the tightness of $U^{n1*}$ and $U^{n2*}$, that is, the tightness of $U^{nijk}$, we deduce that $nM^{nijk}$ has $(\star)$. From (41), the sequence $nN^{nijj}$ is tight and satisfies $(\star)$. Then (c) follows from Theorem 4.2. □

**5. Continuous semimartingales.** Now we consider the case that the driving process $Y$ is a continuous semimartingale. We assume that $Y = H + A$, where $A$ is a continuous adapted process of finite variation and $H$ is a continuous local martingale, both being $d$-dimensional and null at 0. Let $C_t = (C_t^{ij})_{1 \le i,j \le d}$ be the quadratic variation process of $Y$, that is, $C_t^{ij} = [Y^i, Y^j]_t = [H^i, H^j]_t$.

We do not have the necessary and sufficient conditions for the error process $U^n$ to converge at the rate $1/n$. Partial results are available when $A_t$ has the form

$$A_t^i = \int_0^t a_s^i \, ds \qquad \text{with } \int_0^1 |a_s^i|^2 \, ds < \infty \text{ a.s.} \tag{52}$$

THEOREM 5.1. *Assuming that* $Y$ *is a continuous semimartingale, for the following statements, we have that* (a) *and* (52) *imply* (b), *and* (b) *implies* (a).

(a) *There is a* $d \times d$ *nonnegative matrix-valued predictable process* $c$ *such that*

$$C_t = \int_0^t c_s \, ds \qquad \textit{with } \int_0^1 \|c_s\|^3 \, ds < \infty \textit{ a.s.} \tag{53}$$

(b) *The sequence* $(nM^{np}, nN^{np})_{1\le p\le m}$, *in* (6) *and* (7), *has* $(\star)$ *and is tight.*

*And furthermore, under* (a) *and* (52), *if* $H_t = \int_0^t \sigma_s\,dW_s$, *where* $c_t = \sigma_t\sigma_t^\tau$ *as in Theorem* 4.2, *then the sequence* $(nM^{np}, nN^{np})_{1\le p\le d}$ *stably converges in law to the processes* $(M^p, \bar N^p)_{1\le p\le d}$, *where* $\bar N_t^p = N_t^p + 1/2\int_0^t c_s a_s^p\,ds$, *and* $M, N$ *are given in* (50) *and* (51). *Moreover, we have*

$$(Y, nM^{np}, nN^{np})_{1\le p\le d} \text{ has } (\star) \text{ and}$$
$$(Y, nM^{np}, nN^{np})_{1\le p\le d} \Longrightarrow (Y, M^p, \bar N^p)_{1\le p\le d}.$$

Remark 5.1. The limits of the quadratic variation $(nM^n, nN^n, W)$ are listed in Lemma 7.6.

Proof of Theorem 5.1. By the integration by parts formula, we have

$$(Y_s^{(n)})(Y_s^{(n)})^\tau = \int_{n(s)}^s Y_r^{(n)}(\,dY_r)^\tau + \int_{n(s)}^s dY_r(Y_r^{(n)})^\tau + C_s^{(n)}.$$

According to the notation of $M^{np}$ and $N^{np}$ in (6) and (7), we have

$$N^{np} = M^{np} + (M^{np})^\tau + \int_0^t C_s^{(n)}\,dY_s^p.$$

By Theorem 2.3(a) in [2], under (b), we see that the sequence $n\int_0^t C_s^{(n)}\,dY_s^p$ is also tight and has $(\star)$ for any $p = 1, 2, \dots, d$. Particularly, the sequence $n\int_0^t C_s^{(n)pp}\,dY_s^p$ is tight and has $(\star)$. Then its quadratic variation processes, $n^2\int_0^t (C_s^{pp} - C_{n(s)}^{pp})^2\,dC_s^{pp} = n^2 N_t^n(C^{pp})$, are tight too. Theorem 3.2 implies that $\sup_n n^2 N_1^n(C^{pp}) < \infty$. By the same reason as in the proof of (b) $\Rightarrow$ (a) in Theorem 4.2, we can conclude (a).

Now we assume that (a) and (52) hold. By the usual stopping time arguments, we can also assume that $\int_0^1 \|a_s\|^2\,ds$ and $\int_0^1 \|c_s\|^3$ are bounded by a constant $\alpha$.

First, $nN^{np}$ can be decomposed as follows:

$$nN^{np} = F_1^{np} + F_2^{np} + F_3^{np} + F_4^{np} + (F_3^{np} + F_4^{np})^\tau + F_5^{np},$$

where $F_1^{np} = n\int_0^t H_s^{(n)}(H_s^{(n)})^\tau\,dH_s^p$, $F_2^{np} = n\int_0^t H_s^{(n)}(H_s^{(n)})^\tau\,dA_s^p$, $F_3^{np} = n\times \int_0^t H_s^{(n)}(A_s^{(n)})^\tau\,dH_s^p$, $F_4^{np} = n\int_0^t H_s^{(n)}(A_s^{(n)})^\tau\,dA_s^p$, $F_5^{np} = n\int_0^t A_s^{(n)}(A_s^{(n)})^\tau\,dY_s^p$. Since, for any $i$, $n(A_s^{i(n)})^2 = n(A_s^i - A_{n(s)}^i)^2 \le \int_{n(s)}^{n(s)+1/n}(a_s^i)^2\,ds \le \alpha$, by (52) and the Cauchy–Schwarz inequality, $n(A_s^{i(n)})^2$ goes to zero uniformly for $s\in[0,1]$ as $n\to 0$. Thus, $F_5^{np}$ has $(\star)$ and weakly converges to 0. By (81) and (82), $F_3^{np}$ and $F_4^{np}$ have $(\star)$ and weakly converge to 0. By (80) and Theorem 4.2, $F_1^{np} + F_2^{np}$ has $(\star)$ and weakly converges to $\bar N_t^p = N_t^p + \frac12\int_0^t c_s a_s^p\,ds$. Putting them together, $nN^{np}$ has $(\star)$ and weakly converges to $\bar N_t^p$.

Second, $nM^{np}$ can be decomposed similarly as follows:

$$M^{np} = \int_0^t G_1^n\, dH_s^p + \int_0^t G_1^n\, dA_s^p + \int_0^t G_2^n\, dY_s^p + \int_0^t G_3^n\, dY_s^p + \int_0^t G_4^n\, dY_s^p,$$

where $G_1^n = \int_{n(s)}^s H_r^{(n)}\, dH_r^\tau$, $G_2^n = \int_{n(s)}^s H_r^{(n)}\, dA_r^\tau$, $G_3^n = \int_{n(s)}^s A_r^{(n)}\, dH_r^\tau$, $G_4^n = \int_{n(s)}^s A_r^{(n)}\, dA_r^\tau$. When $H$ and $A$ are a one-dimensional process, from (84) and (85), we have $n\int_0^t G_3^n\, dY_s^p \xrightarrow{L^2} 0$; from (81) and (82), we have $n\int_0^t H_s^{(n)} A_s^{(n)}\, dY_s^p \xrightarrow{L^2} 0$; since $G_3^n + G_2^n = H_s^{(n)} A_s^{(n)}$, we have $n\int_0^t G_2^n\, dY_s^p \xrightarrow{L^2} 0$. Therefore, $n\int_0^t G_2^n\, dY_s^p + n\int_0^t G_3^n\, dY_s^p$ has $(\star)$ and weakly converges to 0. Since $n(A_s^{i(n)})^2$ goes to zero uniformly for $s \in [0,1]$ as $n \to \infty$ as shown above, so does $nG_4^n$. Thus, $\int_0^t G_4^n\, dY_s^p$ has $(\star)$ and weakly converges to 0. By (83), $n\int_0^t G_1^n\, dA_s^p$ has $(\star)$ and weakly converges to 0. By Theorem 4.2, we have $n\int_0^t G_1^n\, dY_s^p$ has $(\star)$ and weakly converges to $M_t^p$. So, $nM^{np}$ has $(\star)$ and weakly converges to $M_t^p$. □

Combining Theorems 5.1 and 2.2, we get the rate of convergence $1/n$ for $U^n = X^n - X$, when $Y$ is a continuous semimartingale.

THEOREM 5.2. *Assume that $Y = H + A$ is a continuous semimartingale in $\mathbb{R}^d$ null at zero. If (52) and (53) hold, then for any starting point $x_0$ and any $C^2$ functions $f$ with at most linear growth, the sequence $(Y, nU^n)$ converges weakly to a limit $(Y, U)$, the process $U = (U^1, \dots, U^q)^\tau$ is the solution of the linear equation (8), with $M$ given in (50) and $N$ replaced by $\bar{N}$, which is given in Theorem 5.1.*

**6. Examples.** For the Itô-type SDE in (2), $X_t = x_0 + \int_0^t a(X_s)\, dW_s + \int_0^t b(X_s)\, ds$. If $a, b$ are $C^2$ functions, then we can use the Milstein scheme $X_t^n$ defined in (4) to approximate the solution $X_t$,

$$X_t^n = X_{n(t)}^n + aW_t^{(n)} + bt^{(n)} + \tfrac{1}{2}aa'((W_t^{(n)})^2 - t^{(n)}) \\ + \tfrac{1}{2}bb'(t^{(n)})^2 + \int_{n(t)}^t ab' s^{(n)}\, dW_s + \int_{n(t)}^t ab' W_s^{(n)}\, ds, \tag{54}$$

where $W_t^{(n)} = W_t - W_{n(t)}$, $t^{(n)} = t - n(t)$, and function $a, a', b, b'$ taking values at $X_{n(t)}^n$. Compared with the Milstein scheme defined in (3), this scheme has two more terms, $\int_{n(t)}^t ab' s^{(n)}\, dW_s$ and $\int_{n(t)}^t ab' W_s^{(n)}\, ds$.

In order to apply Theorem 5.1 to find the asymptotic error $U$ of $nU^n = n(X^n - X)$ for the scheme in (54), we use the following notation: $f =$

$(a,b)$, $Y_t=(W_t,t)^\tau$, $H_t=(W_t,0)^\tau$, $A_t=(0,t)^\tau$. Let $F_s=Y_s^{(n)}(Y_s^{(n)})^\tau$, $G_s=\int_{n(s)}^s Y_r^{(n)}d(Y_r)^\tau$, and

$$F_s=\begin{pmatrix} W_s^{(n)}W_s^{(n)} & W_s^{(n)}s^{(n)} \\ s^{(n)}W_s^{(n)} & s^{(n)}s^{(n)} \end{pmatrix},\qquad G_s=\begin{pmatrix} \int_{n(s)}^s W_r^{(n)}\,dW_r & \int_{n(s)}^s W_r^{(n)}\,dr \\ \int_{n(s)}^s r^{(n)}\,dW_r & \int_{n(s)}^s r^{(n)}\,dr \end{pmatrix}$$

and

$$N_t^{n1}=\int_0^t F_s\,dW_s,\qquad N_t^{n2}=\int_0^t F_s\,ds,$$
$$M_t^{n1}=\int_0^t G_s\,dW_s,\qquad M_t^{n2}=\int_0^t G_s\,ds.$$

By Theorem 5.1, $(nN^{n111},nN^{n211},nM^{n111})\Longrightarrow(N^{111},N^{211},M^{111})$, where

$$N_t^{211}=\frac{t}{2},\qquad M^{111}=\frac{1}{\sqrt{6}}B_t^1,\qquad N_t^{111}=\frac{2}{\sqrt{6}}B_t^1+\frac{1}{2\sqrt{3}}B_t^2+\frac{1}{2}W_t,$$

and $B^1,B^2,W$ are three independent standard Brownian motions, defined on an extension of the probability space on which $W$ is defined. We also have that the limits of all other entries of $nN^{n1},nN^{n2},nM^{n1},nM^{n2}$ are zeros. Therefore, the limit process $U$ of $n(X^n-X)$ satisfies the following linear SDE:

$$\begin{aligned} U_t=&\int_0^t U_s a'(X_s)\,dW_s+U_s b'(X_s)\,ds-\frac{1}{4}\int_0^t a^2(X_s)b''(X_s)\,ds\\ &-\int_0^t a(X_s)(a'(X_s))^2\frac{dB_s^1}{\sqrt{6}}-\int_0^t a^2(X_s)a''(X_s)\left(\frac{dB_s^1}{\sqrt{6}}+\frac{dB_s^2}{\sqrt{12}}+\frac{dW_s}{4}\right).\end{aligned}$$

If $a(\cdot)$ is a $C^2$ function and $b(\cdot)$ is a $C^1$ function, we can use the following simpler Milstein scheme $X^n$ rather than (54):

$$X_t^n=x_0+\int_0^t b(X_{n(s)}^n)\,ds+\int_0^t a(X_{n(s)}^n)\,dW_s+\int_0^t a(X_{n(s)}^n)a'(X_{n(s)}^n)W_s^{(n)}\,dW_s$$

to achieve the same rate of convergence, that is, $n(X^n-X)\Rightarrow U$, and $U$ satisfies

$$\begin{aligned} U_t=&\int_0^t U_s b'(X_s)\,ds+U_s a'(X_s)\,dW_s-\frac{1}{2}\int_0^t b(X_s)b'(X_s)\,ds\\ &-\frac{1}{2}\int_0^t C_0(X_s)\,dW_s-\frac{1}{\sqrt{12}}\int_0^t C_1(X_s)\,dB_s^1-\frac{1}{\sqrt{6}}\int_0^t C_2(X_s)\,dB_s^2,\end{aligned}\tag{55}$$

where $W,B^1,B^2$ are three independent Brownian motions and $C_0,C_1,C_2$ are three functions dependent on $a$ and $b$.

**7. Some lemmas.** In the next two lemmas we assume that $X, Y$ and $Z$ are three continuous processes with finite variations, which satisfy

$$X_t = \int_0^t x_s\,ds, \qquad Y_t = \int_0^t y_s\,ds, \qquad Z_t = \int_0^t z_s\,ds, \qquad x, y, z \in L^3([0,1]), \tag{56}$$

where $L^p([0,1]) = \{x \| \|x\|_p = (\int_0^1 |x_s|^p\,dx)^{1/p} < \infty\}$.

LEMMA 7.1. *Let $X_s^{(n)} = X_s - X_{n(s)}$, then for any $0 \le t \le 1$,*

$$n^2 \left| \int_0^t X_s^{(n)} Y_s^{(n)}\,dZ_s \right| \le \|x\|_3 \|y\|_3 \|z\|_3, \tag{57}$$

$$n^2 \left| \int_0^t \int_{n(s)}^s X_r^{(n)}\,dY_r\,dZ_s \right| \le \|x\|_3 \|y\|_3 \|z\|_3. \tag{58}$$

PROOF. By the Hölder inequality, $|X_s^{(n)}|^3 \le (s - n(s))^2 \int_{n(s)}^s |x_r|^3\,dr$. Then

$$\|X^{(n)}\|_3^3 \le \int_0^1 (s - n(s))^2 \int_{n(s)}^s |x_r|^3\,dr\,ds \le \frac{1}{n^3} \|x\|_3^3.$$

It follows that $n^2 |\int_0^t X_s^{(n)} Y_s^{(n)}\,dZ_s| \le n^2 \|X^{(n)}\|_3 \|Y^{(n)}\|_3 \|z\|_3 \le \|x\|_3 \|y\|_3 \|z\|_3$, which proves (57). Since

$$\left| \int_{n(s)}^s X_r^{(n)}\,dY_r \right| \le \int_{n(s)}^s \int_{n(r)}^r |x_u|\,du\,|y_s|\,ds \le \int_{n(s)}^s |x_u|\,du \int_{n(s)}^s |y_r|\,dr, \tag{59}$$

(58) follows. □

LEMMA 7.2. *As $n \to \infty$, the following uniform convergences in $t$ holds:*

$$n^2 \int_0^t X_s^{(n)} X_s^{(n)}\,dZ_s \longrightarrow \tfrac{1}{3} \int_0^t x_s x_s z_s\,ds, \tag{60}$$

$$n^2 \int_0^t X_s^{(n)} Y_s^{(n)}\,dZ_s \longrightarrow \tfrac{1}{3} \int_0^t x_s y_s z_s\,ds, \tag{61}$$

$$n^2 \int_0^t \int_{n(s)}^s X_r^{(n)} dY_r\,dZ_s \longrightarrow \tfrac{1}{6} \int_0^t x_s y_s z_s\,ds. \tag{62}$$

REMARK 7.1. When $x, z \in L^2([0,1])$, from Theorem 4.1 in [2], we have

$$n \int_0^t X_s^{(n)}\,dZ_s \longrightarrow \tfrac{1}{2} \int_0^t x_s z_s\,ds \qquad \text{uniformly in } t \text{ as } n \to \infty. \tag{63}$$

When $\|x\|_2$ and $\|z\|_2$ are bounded, we have $L^2$ convergence in (63).

REMARK 7.2. By Cauchy's inequality, we have the following, which will be used in the proof of Lemma 7.7:

$$n\int_0^1 |X_s^{(n)}|\,ds \le \int_0^1 |x_s|\,ds \quad \text{and} \quad n^2\int_0^1 |X_s^{(n)}|^2\,ds \le \int_0^1 |x_s|^2\,ds. \tag{64}$$

PROOF OF LEMMA 7.2. To prove (60), first we assume that $x_s \in C^1$, therefore, $x_s$ and $dx_s/ds$ are bounded by a constant $\alpha$ in $[0,1]$. By the mean value theorem, there exists an $\bar{s}$ between $n(s)$ and $s$ such that $X_s^{(n)} = x_{\bar{s}}(s-n(s))$. Since $n^2\int_0^t (s-n(s))^2\,ds \to t/3$, and

$$n^2\int_0^t (x_s)^2(s-n(s))^2 z_s\,ds \longrightarrow \tfrac{1}{3}\int_0^t (x_s)^2 z_s\,ds,$$

we need only to show that $n^2\int_0^t((x_s)^2-(x_{\bar{s}})^2)(s-n(s))^2 z_s\,ds \to 0$. This is true because $|(x_s)^2-(x_{\bar{s}})^2| \le 2\alpha^2/n$ and $n^2\int_0^t (s-n(s))^2 z_s\,ds \to \frac{1}{3}\int_0^t z_s\,ds$. Second, if $x_s \notin C^1$, since $x_s \in L^3$, for any $\varepsilon > 0$, there exist a $\bar{x}_s \in C^1$ such that $\|x-\bar{x}\|_3 \le \varepsilon$. By (57),

$$\begin{aligned}
&\left|n^2\int_0^t \left(\int_{n(s)}^s x_r\,dr\right)^2 dZ_s - n^2\int_0^t \left(\int_{n(s)}^s \bar{x}_r\,dr\right)^2 dZ_s\right| \\
&\quad = n^2\left|\int_0^t \int_{n(s)}^s (x_r+\bar{x}_r)\,dr \int_{n(s)}^s (x_r-\bar{x}_r)\,dr\,dZ_s\right| \\
&\quad \le \|x+\bar{x}\|_3\|x-\bar{x}\|_3\|z\|_3.
\end{aligned}$$

Now we can claim that (60) holds for all $x \in L^3$. Since $2ab = (a+b)^2 - a^2 - b^2$, (61) follows by (60). (62) can be proved by the same arguments as in the proof of (60). Here we just give two key steps. It is easy to see that (62) holds when $x_s, y_s \in C^1$, since there exist $\bar{s}$ and $\bar{r}$ between $n(s)$ and $s$ such that $\int_{n(s)}^s X_r^{(n)}\,dY_r = x_{\bar{s}}y_{\bar{r}}\int_{n(s)}^s dr\int_{n(r)}^r du = x_{\bar{s}}y_{\bar{r}}(s-n(s))^2/2$. And for $x, y, \bar{x}, \bar{y}$, by (58),

$$\begin{aligned}
&\left|n^2\int_0^t\int_{n(s)}^s\int_{n(r)}^r x_u\,du\,y_r\,dr\,dZ_s - n^2\int_0^t\int_{n(s)}^s\int_{n(r)}^r \bar{x}_u\,du\,\bar{y}_r\,dr\,dZ_s\right| \\
&\quad \le \|x\|_3\|y-\bar{y}\|_3\|z\|_3 + \|\bar{y}\|_3\|x-\bar{x}\|_3\|z\|_3.
\end{aligned}$$

If $x_s, y_s$ and $z_s$ are all nonnegative, both sides of (60), (61) and (62) are all continuous and nondecreasing, these convergences are uniform in $t$ over $[0,1]$. This uniformity can be easily generalized for all $X, Y$ and $Z$ which satisfy (56). $\square$

In the proof of the next two lemmas, let $C_p$ be the constant in Burkholder's inequality, that is, $\mathbb{E}|Y_t|^p \le C_p\mathbb{E}[Y,Y]_t^{p/2}$ for any continuous local martingale $Y$.

LEMMA 7.3. *Let $W, B, U$ and $V$ be four independent standard Brownian motions and $n(s) = [ns]/n$. If we write $W_s^{(n)} = W_s - W_{n(s)}$, then for $0 \le t \le 1$,*

(a) $n^2 \int_0^t W_s^{(n)} W_s^{(n)} W_s^{(n)} W_s^{(n)}\, ds \xrightarrow{L^2} t$,

(b) $n^2 \int_0^t W_s^{(n)} W_s^{(n)} B_s^{(n)} B_s^{(n)}\, ds \xrightarrow{L^2} t/3$,

(c) $n^2 \int_0^t W_s^{(n)} W_s^{(n)} W_s^{(n)} B_s^{(n)}\, ds \xrightarrow{L^2} 0$,

(d) $n^2 \int_0^t W_s^{(n)} W_s^{(n)} U_s^{(n)} V_s^{(n)}\, ds \xrightarrow{L^2} 0$,

(e) $n^2 \int_0^t W_s^{(n)} B_s^{(n)} U_s^{(n)} V_s^{(n)}\, ds \xrightarrow{L^2} 0$.

PROOF. We prove the lemma by the simple fact that a sequence of random variables $\xi_n \xrightarrow{L^2}$ the limit of $\mathbb{E}[\xi_n]$, if the variance of $\xi_n$, $\mathrm{Var}(\xi_n)$, converges to zero as $n \to \infty$. Now for any four standard Brownian motions $W^1, W^2, W^3, W^4$ (not necessary to be independent), let $A_s^{(n)} = W_s^{1(n)} W_s^{2(n)} W_s^{3(n)} W_s^{4(n)}$; by Cauchy's inequality, we have $\mathbb{E}(A_s^{(n)})^2 \le \mathbb{E}(W_s^{(n)})^8 \le C_8[s - n(s)]^4$. Therefore,

$$\begin{aligned}\mathrm{Var}\biggl(n^2 \int_0^t A_s^{(n)}\, ds\biggr) &\le n^5 \mathbb{E}\biggl(\int_0^{1/n} A_s^{(n)}\, ds\biggr)^2 \\ &\le n^4 \int_0^{1/n} \mathbb{E}(A_s^{(n)})^2\, ds \\ &\le C_8 n^4 \int_0^{1/n} s^4\, ds \le C_8/n,\end{aligned}$$

which implies that the variances of the left-hand sides of (a), (b), (c), (d) and (e) converge to zero as $n \to \infty$. Next we calculate the expectations of the left-hand sides. Since the expectations of the left-hand sides of (c), (d) and (e) are 0 and

$$\begin{aligned}\mathbb{E} n^2 \int_0^t (W_s^{(n)})^4\, ds &= n^2 \int_0^t 3(s - [ns]/n)^2\, ds \\ &= [nt]/n + (nt - [nt])/n^3 \longrightarrow t, \\ \mathbb{E} n^2 \int_0^t (W_s^{(n)} B_s^{(n)})^2\, ds &= n^2 \int_0^t (s - [ns]/n)^2\, ds \longrightarrow t/3,\end{aligned}$$

we can conclude all of the $L^2$ convergences. □

LEMMA 7.4. *Suppose that for any pair from four standard Brownian motions $W, B, U$ and $V$, they are either identical or independent, then for* $0 \le t \le 1$,

$$\text{(a)}\quad n^2 \int_0^t \left\{ \int_{n(s)}^s W_r^{(n)}\, dB_r \int_{n(s)}^s U_r^{(n)}\, dV_r \right\} ds \xrightarrow{L^2} \begin{cases} t/6, & \textit{if } W = U \textit{ and } B = V, \\ 0, & \textit{otherwise;} \end{cases}$$

$$\text{(b)}\quad n^2 \int_0^t \left\{ W_s^{(n)} B_s^{(n)} \int_{n(s)}^s U_r^{(n)}\, dV_r \right\} ds \xrightarrow{L^2} \begin{cases} t/3, & \textit{if } W = B = U = V, \\ t/6, & \textit{if } W = U \ne B = V, \\ t/6, & \textit{if } B = U \ne W = V, \\ 0, & \textit{otherwise.} \end{cases}$$

PROOF. We use the same arguments as in the proof of Lemma 7.3 to prove the $L^2$ convergence. Let

$$X_s = \int_{n(s)}^s W_r^{(n)}\, dB_r \quad \text{and} \quad Y_s = \int_{n(s)}^s U_r^{(n)}\, dV_r.$$

First, for the variance of the left-hand side of (a), we have

$$\begin{aligned} \mathbb{E}(X_s)^4 &\le C_4 \mathbb{E}\left( \int_{n(s)}^s [W_r^{(n)}]^2\, dr \right)^2 \\ &\le C_4 (s - n(s)) \mathbb{E} \int_{n(s)}^s [W_r^{(n)}]^4\, dr \\ &= C_4 (s - n(s))^4. \end{aligned}$$

Similarly, $\mathbb{E}(Y_s)^4 \le C_4(s - n(s))^4$. By the Cauchy inequality,

$$\begin{aligned} \operatorname{Var}\left( n^2 \int_0^t X_s Y_s\, ds \right) &\le n^5 \mathbb{E}\left( \int_0^{1/n} X_s Y_s\, ds \right)^2 \le n^4 \mathbb{E} \int_0^{1/n} (X_s Y_s)^2\, ds \\ &\le n^4 \int_0^{1/n} (\mathbb{E}(X_s)^4 \mathbb{E}(Y_s)^4)^{1/2}\, ds \\ &\le n^4 \int_0^{1/n} C_4 (s - n(s))^4\, ds = C_4/(5n) \to 0, \end{aligned}$$

as $n \to \infty$. For the expectation of the left-hand side of (a), we use the integration by parts formula (see page 60 of [10]). If $B$ and $V$ are independent, then $\mathbb{E}(X_s Y_s) = \mathbb{E}[X, Y]_s = 0$. If $B = V$, then $\mathbb{E}(X_s Y_s) = \mathbb{E}[X, Y]_s = \mathbb{E} \int_{n(s)}^s W_r^{(n)} U_r^{(n)}\, dr$. Now if $W$ and $U$ are independent, then $\mathbb{E}(W_r^{(n)} U_r^{(n)}) = 0$; if $W = U$, then $\mathbb{E}(W_r^{(n)} U_r^{(n)}) = r - n(r)$. Therefore, if $B = V$ and $W = U$,

$$\lim_{n\to\infty} \mathbb{E} n^2 \int_0^t X_s Y_s\, ds = \lim_{n\to\infty} n^2 \int_0^t \int_{n(s)}^s (r - n(r))\, dr\, ds$$

$$= \lim_{n\to\infty} n^2/2 \int_0^t (s-n(s))^2\, ds = t/6;$$

otherwise, it is zero. We finish the proof of (a). Second, we claim that the variance of the left-hand side in (b) goes to zero as $n \to \infty$. This is because

$$\begin{aligned}
&\operatorname{Var}\biggl(n^2 \int_0^t W_s^{(n)} B_s^{(n)} \int_{n(s)}^s U_r^{(n)}\, dV_r\, ds\biggr) \\
&\quad \le n^5 \mathbb{E}\biggl(\int_0^{1/n} W_s^{(n)} B_s^{(n)} Y_s\, ds\biggr)^2 \le n^4 \mathbb{E} \int_0^{1/n} (W_s^{(n)} B_s^{(n)} Y_s)^2\, ds \\
&\quad \le n^4 \int_0^{1/n} \sqrt{\mathbb{E}(W_s^{(n)} B_s^{(n)})^4 \mathbb{E}(Y_s)^4}\, ds \le n^4 \int_0^{1/n} \sqrt{\mathbb{E}(W_s^{(n)})^8 \mathbb{E}(Y_s)^4}\, ds \\
&\quad \le \sqrt{C_4 C_8} n^4 \int_0^{1/n} (s-n(s))^4\, ds \le \sqrt{C_4 C_8}/(5n).
\end{aligned}$$

Next, we prove the convergence of the expectation of the left-hand side of (b) to the right-hand side of (b). For the first case that $W = B = U = V$, since $(W_s^{(n)})^2 = 2\int_{n(s)}^s W_r^{(n)}\, dW_r + s - n(s)$,

$$\begin{aligned}
\mathbb{E}\biggl\{n^2 \int_0^t (W_s^{(n)})^2 \int_{n(s)}^s W_r^{(n)}\, dW_r\, ds\biggr\} &= 2n^2 \int_0^t \mathbb{E}\biggl\{\int_{n(s)}^s W_r^{(n)}\, dW_r\biggr\}^2 ds \\
&= n^2 \int_0^t (s-n(s))^2\, ds \to t/3 \\
&\qquad \text{as } n \to \infty.
\end{aligned}$$

For the second case that $W = U \ne B = V$, since $W$ and $B$ are independent, $W_s^{(n)} B_s^{(n)} = \int_{n(s)}^s W_r^{(n)}\, dB_r + \int_{n(s)}^s B_r^{(n)}\, dW_r$; taking expectations of the squares of both sides, we have $\mathbb{E}\{\int_{n(s)}^s W_r^{(n)}\, dB_r \int_{n(s)}^s B_r^{(n)}\, dW_r\} = 0$. Therefore,

$$\begin{aligned}
\mathbb{E}\biggl\{n^2 \int_0^t W_s^{(n)} B_s^{(n)} \int_{n(s)}^s W_r^{(n)}\, dB_r\, ds\biggr\} &= n^2 \int_0^t \mathbb{E}\biggl(\int_{n(s)}^s W_r^{(n)}\, dB_r\biggr)^2 ds \\
&= n^2/2 \int_0^t (s-n(s))^2\, ds \to t/6 \\
&\qquad \text{as } n \to \infty.
\end{aligned}$$

The third case is equivalent to the second. It is easy to see that the expectations of both sides in the forth case are zeros. □

Now, let $Y$ be the continuous local martingale defined as the beginning of Section 4, that is, $Y_t = Y_t(\sigma) = \int_0^t \sigma_s\, dW_s$. In order to prove Theorem 4.2 in the general case of $\sigma$, we define another continuous local martingale $\bar{Y}_t$, on the probability space where $Y$ and $W$ live on, as follows: $\bar{Y}_t = Y_t(\bar{\sigma}) = \int_0^t \bar{\sigma}_s \, dW_s$, where $\bar{\sigma}$ is a $d \times m$ matrix-valued process. We write

$$N_t^{np}(\sigma) = \int_0^t Y_s^{(n)}(\sigma)(Y_s^{(n)}(\sigma))^\tau \, dY_s^p(\sigma), \tag{65}$$

$$M_t^{np}(\sigma) = \int_0^t \int_{n(s)}^s Y_r^{(n)}(\sigma) \, dY_r{}^\tau(\sigma) \, dY_s^p(\sigma),$$

$$\eta = \eta(\sigma, \bar{\sigma}) = \left(\mathbb{E}\int_0^1 \|\sigma_s - \bar{\sigma}_s\|^6 \, ds\right)^{1/6}, \tag{66}$$

where $\|\sigma\| = (\sum_{ij} [\sigma_s^{ij}]^2)^{1/2}$. $N_t^{np}(\bar{\sigma})$ and $M_t^{np}(\bar{\sigma})$ are defined similarly for $Y_t(\bar{\sigma})$.

LEMMA 7.5. *If $\int_0^1 \|\sigma\|^6 \, ds$ and $\int_0^1 \|\bar{\sigma}\|^6 \, ds$ have an upper bound $p$, then there exists a constant $C$, which depends on $p$, such that*

$$n^2 \mathbb{E}\|[N^{np}(\sigma), N^{nq}(\sigma)]_t - [N^{np}(\bar{\sigma}), N^{nq}(\bar{\sigma})]_t\| \le C\eta, \tag{67}$$

$$n^2 \mathbb{E}\|[M^{np}(\sigma), M^{nq}(\sigma)]_t - [M^{np}(\bar{\sigma}), M^{nq}(\bar{\sigma})]_t\| \le C\eta, \tag{68}$$

$$n^2 \mathbb{E}\|[N^{np}(\sigma), M^{nq}(\sigma)]_t - [N^{np}(\bar{\sigma}), M^{nq}(\bar{\sigma})]_t\| \le C\eta, \tag{69}$$

$$n^2 \mathbb{E}\|[N^{np}(\sigma), W]_t - [N^{np}(\bar{\sigma}), W]_t\| \le C\eta, \tag{70}$$

$$n^2 \mathbb{E}\|[M^{np}(\sigma), W]_t - [M^{np}(\bar{\sigma}), W]_t\| \le C\eta. \tag{71}$$

PROOF. Let $C$ be the constant which changes from line to line. By the Burkholder and the Hölder inequalities, for an integer $k$,

$$\mathbb{E}(\|Y_s^{(n)}\|)^{2k} \le C(s - n(s))^{k-1} \mathbb{E}\int_{n(s)}^s \|\sigma_r\|^{2k} \, dr, \tag{72}$$

thus, $\int_0^t \mathbb{E}(\|Y_s^{(n)}\|)^{2k} \, ds \le C/n^k$ and also $\int_0^t \mathbb{E}(\|\bar{Y}_s^{(n)}\|)^{2k} \, ds \le C/n^k$. Similarly,

$$\int_0^t \mathbb{E}(\|Y_s^{(n)} - \bar{Y}_s^{(n)}\|)^{2k} \, ds \le C/n^k \mathbb{E}\int_0^t \|\sigma_r - \bar{\sigma}_r\|^{2k} \, dr.$$

Let $A \otimes B = (a_{ij} B)$ be the Kronecker product between two matrices $A = (a_{ij})$ and $B = (b_{nm})$. For (67), we write $A = Y_s^{(n)}(Y_s^{(n)})^\tau$ and $B = \bar{Y}_s^{(n)}(\bar{Y}_s^{(n)})^\tau$, then by the Hölder inequality,

$$\begin{aligned} n^2 \mathbb{E}\|[N^{np}(\sigma), & N^{nq}(\sigma)]_t - [N^{np}(\bar{\sigma}), N^{nq}(\bar{\sigma})]_t\| \\ &= n^2 \mathbb{E}\left\|\int_0^t A \otimes A c^{pq}(\sigma) - B \otimes B c^{pq}(\bar{\sigma}) \, ds\right\| \\ &\le n^2 \mathbb{E}\int_0^t \|Y_s^{(n)}\|^4 \|\sigma - \bar{\sigma}\|(\|\sigma\| + \|\bar{\sigma}\|) \end{aligned}$$

$$+ 6\|Y_s^{(n)} - \bar{Y}_s^{(n)}\|(\|Y_s^{(n)}\|^3 + \|\bar{Y}_s^{(n)}\|^3)\|\bar{\sigma}\|^2\, ds$$

$$\leq Cn^2\biggl(\int_0^t \mathbb{E}\|Y_s^{(n)} - \bar{Y}_s^{(n)}\|^6\, ds\biggr)^{1/6}\biggl(\int_0^t \mathbb{E}(\|Y_s^{(n)}\|^3 + \|\bar{Y}_s^{(n)}\|^3)^2\, ds\biggr)^{1/2}$$

$$+ Cn^2\biggl(\int_0^t \mathbb{E}\|Y_s^{(n)}\|^6\, ds\biggr)^{2/3}\eta \leq C\eta.$$

For (68), we write $A = \int_{n(s)}^s Y_r^{(n)}\, d(Y_r)^\tau$ and $B = \int_{n(s)}^s \bar{Y}_r^{(n)}\, d(\bar{Y}_r)^\tau$, again, by the Burkholder and the Hölder inequalities and (72),

$$\begin{aligned}
[\mathbb{E}\|A\|^3]^2 &\leq C\biggl[\mathbb{E}\biggl(\int_{n(s)}^s \|Y_r^{(n)}\|^2\|\sigma_r\|^2\, dr\biggr)^{3/2}\biggr]^2 \\
&\leq C\biggl[\mathbb{E}\biggl(\int_{n(s)}^s \|Y_r^{(n)}\|^3\, dr\biggr)\biggl(\int_{n(s)}^s \|\sigma_r\|^6\, dr\biggr)^{1/2}\biggr]^2 \\
&\leq C\mathbb{E}\biggl(\int_{n(s)}^s \|Y_r^{(n)}\|^3\, dr\biggr)^2\mathbb{E}\biggl(\int_{n(s)}^s \|\sigma_r\|^6\, dr\biggr) \\
&\leq C\biggl((s - n(s))\int_{n(s)}^s [r - n(r)]^2\mathbb{E}\biggl[\int_{n(r)}^r \|\sigma_u\|^6\, du\biggr]dr\biggr) \\
&\qquad \times \biggl(\mathbb{E}\biggl[\int_{n(s)}^s \|\sigma_r\|^6\, dr\biggr]\biggr) \\
&\leq C(s - n(s))^4\biggl(\mathbb{E}\biggl[\int_{n(s)}^s \|\sigma_r\|^6\, dr\biggr]\biggr)^2.
\end{aligned}$$

It follows that

$$\begin{aligned}
\int_0^t \mathbb{E}\|A\|^3\, ds &\leq C\int_0^t (s - n(s))^2\mathbb{E}\biggl[\int_{n(s)}^s \|\sigma_r\|^6\, dr\biggr]\, ds \\
&\leq 1/n^3\mathbb{E}\int_0^t \|\sigma_s\|^6\, ds.
\end{aligned} \tag{73}$$

Similarly,

$$\begin{aligned}
\mathbb{E}\|A - B\|^2 &\leq C\mathbb{E}\int_{n(s)}^s \|\sigma_r Y_r^{(n)} - \bar{\sigma}_r\bar{Y}_r^{(n)}\|^2\, dr \\
&\leq C\mathbb{E}\int_{n(s)}^s \|Y_r^{(n)} - \bar{Y}_r^{(n)}\|^2\|\sigma_r\|^2 + \|\bar{Y}_r^{(n)}\|^2\|\sigma_r - \bar{\sigma}_r\|^2\, dr \\
&\leq C\biggl([s - n(s)]^2\mathbb{E}\int_{n(s)}^s \|\sigma_r - \bar{\sigma}\|^4\, dr\mathbb{E}\int_{n(s)}^s \|\sigma_r\|^4 + \|\bar{\sigma}_r\|^4\, dr\biggr)^{1/2}.
\end{aligned}$$

By the Cauchy inequality,

$$\left(\int_0^t \mathbb{E}\|A-B\|^2\,ds\right)^2 \le C\int_0^t [s-n(s)]^2\mathbb{E}\int_{n(s)}^s \|\sigma_r-\bar\sigma\|^4\,dr\,ds$$
$$\times\int_0^t \mathbb{E}\int_{n(s)}^s \|\sigma_r\|^4+\|\bar\sigma_r\|^4\,dr\,ds \tag{74}$$
$$\le C1/n^4\mathbb{E}\int_0^1 \|\sigma_s-\bar\sigma_s\|^4\,ds.$$

Finally, (68) follows from (73), (74) and

$$n^2\mathbb{E}\|[M^{np}(\sigma),M^{nq}(\sigma)]_t-[M^{np}(\bar\sigma),M^{nq}(\bar\sigma)]_t\|$$
$$= n^2\mathbb{E}\left\|\int_0^t A\otimes Ac^{pq}(\sigma)-B\otimes Bc^{pq}(\bar\sigma)\,ds\right\|$$
$$\le n^2\mathbb{E}\int_0^t \|A\|^2\|\sigma-\bar\sigma\|(\|\sigma\|+\|\bar\sigma\|)+\|A-B\|(\|A\|+\|B\|)\|\bar\sigma\|^2\,ds$$
$$\le n^2C\left(\int_0^t \mathbb{E}\|A\|^3\,ds\right)^{2/3}\eta$$
$$+n^2C\left(\int_0^t \mathbb{E}\|A-B\|^2\,ds\right)^{1/2}\left(\int_0^t \mathbb{E}(\|A\|+\|B\|)^3\,ds\right)^{1/3}.$$

Inequalities (69), (70) and (71) can be proved in the same way as above. □

LEMMA 7.6. *Let $N^{npij}(\sigma)$ and $M^{npij}(\sigma)$ be the $(i,j)$ entries of $N^{np}(\sigma)$ and $M^{np}(\sigma)$, respectively. Let $W^a$ be the $a$th component of the Brownian motion $W$ in $\mathbb{R}^d$. Under the assumption of (49), for all $1\le i,j,k,l,p,q\le d$ and $1\le a\le m$, and $t\in[0,1]$, we have*

$$n^2[N^{npij}(\sigma),N^{nqkl}(\sigma)]_t \xrightarrow{P} \tfrac13\int_0^t (c_s^{ij}c_s^{kl}+c_s^{ik}c_s^{jl}+c_s^{il}c_s^{jk})c_s^{pq}\,ds, \tag{75}$$

$$n^2[N^{npij}(\sigma),M^{nqkl}(\sigma)]_t \xrightarrow{P} \tfrac16\int_0^t (c_s^{ik}c_s^{jl}+c_s^{il}c_s^{jk})c_s^{pq}\,ds, \tag{76}$$

$$n^2[M^{npij}(\sigma),M^{nqkl}(\sigma)]_t \xrightarrow{P} \tfrac16\int_0^t c_s^{ik}c_s^{jl}c_s^{pq}\,ds, \tag{77}$$

$$n[N^{np\,ij}(\sigma),W^a]_t \xrightarrow{P} \tfrac12\int_0^t c_s^{ij}\sigma_s^{pk}\,ds, \tag{78}$$

$$n[M^{npij}(\sigma),W^a]_t \xrightarrow{P} 0. \tag{79}$$

PROOF. Let $T_p=\inf\{t:\int_0^t\|c_s\|^3ds\ge p\}$. Since (49) implies that $P(T_p<1)\to 0$, it is enough to prove the lemma for the processes stopped at $T_p$, which is equivalent to assuming that $\int_0^1 \|c_s\|^3\,ds$ is bounded by a constant $p$, or $\int_0^1 \|\sigma_s\|^6\,ds \le p$. Since $\sigma_s \in L^6([0,1])$, for any $\varepsilon > 0$, there exists a matrix function $\tilde{\sigma}_s \in C^6([0,1])$ such that $\int_0^1 \|\sigma_s - \tilde{\sigma}_s\|^6\,ds \le \varepsilon$, and $\int_0^1 \|\tilde{\sigma}_s\|^6\,ds \le p$. If we denote the right-hand side of (75) as $h(\sigma)$, and denote it as $h(\tilde{\sigma})$ when $c_s$ is replaced by $\tilde{\sigma}_s\tilde{\sigma}_s^\tau$, then $|h(\sigma) - h(\tilde{\sigma})| \le \int_0^1 \|\sigma_s - \tilde{\sigma}_s\|^6\,ds$. We have the same results for the right-hand sides of (76), (77), (78) and (79). Therefore, by Lemma 7.5, it suffices to prove the lemma for $\sigma_s \in C^6([0,1])$.

Now we define $\bar{\sigma}_s = \sigma_{n(s)}$, where $n(s) = [ns]/n$. Since $\int_0^1 \|\sigma_s - \bar{\sigma}_s\|^6\,ds \le \kappa/n$ for some constant $\kappa$, by Lemma 7.5 again, we need only to show that

$$[N^{npij}(\bar{\sigma}), N^{nqkl}(\bar{\sigma})]_t \xrightarrow{P} \frac{1}{3}\int_0^t (c_s^{ij}c_s^{kl} + c_s^{ik}c_s^{jl} + c_s^{il}c_s^{jk})c_s^{pq}\,ds,$$

and the corresponding convergences in (76), (77), (78) and (79) where $\sigma$ is replaced by $\bar{\sigma}$. By the monotone class theorem and Lemmas 7.3 and 7.4,

(a) $$n^2[N^{npij}(\bar{\sigma}), N^{nqkl}(\bar{\sigma})]_t = \sum_{i'j'k'l'} n^2 \int_0^t \bar{\sigma}_s^{ii'}\bar{\sigma}_s^{jj'}\bar{\sigma}_s^{kk'}\bar{\sigma}_s^{ll'} c_s^{pq}(\bar{\sigma}) W_s^{(n)i'} W_s^{(n)j'} W_s^{(n)k'} W_s^{(n)l'}\,ds \xrightarrow{P} \frac{1}{3}\int_0^t (c_s^{ij}c_s^{kl} + c_s^{ik}c_s^{jl} + c_s^{il}c_s^{jk})c_s^{pq}\,ds,$$

(b) $$n^2[N^{npij}(\bar{\sigma}), M^{nqkl}(\bar{\sigma})]_t = \sum_{i'j'k'l'} n^2 \int_0^t \bar{\sigma}_s^{ii'}\bar{\sigma}_s^{jj'}\bar{\sigma}_s^{kk'}\bar{\sigma}_s^{ll'} c_s^{pq}(\bar{\sigma}) W_s^{(n)i'} W_s^{(n)j'} \int_{n(s)}^s W_r^{(n)k'}\,dW_r^{l'}\,ds \xrightarrow{P} \frac{1}{6}\int_0^t (c_s^{ik}c_s^{jl} + c_s^{il}c_s^{jk})c_s^{pq}\,ds,$$

(c) $$n^2[M^{npij}(\bar{\sigma}), M^{nqkl}(\bar{\sigma})]_t = \sum_{i'j'k'l'} n^2 \int_0^t \bar{\sigma}_s^{ii'}\bar{\sigma}_s^{jj'}\bar{\sigma}_s^{kk'}\bar{\sigma}_s^{ll'} c_s^{pq}(\bar{\sigma}) \times \int_{n(s)}^s W_r^{(n)i'}\,dW_r^{j'} \int_{n(s)}^s W_r^{(n)k'}\,dW_r^{l'}\,ds \xrightarrow{P} \frac{1}{6}\int_0^t c_s^{ik}c_s^{jl}c_s^{pq}\,ds,$$

(d) $$n[N^{npij}(\bar{\sigma}), W^k]_t = \sum_{i'j'} n \int_0^t \bar{\sigma}_s^{ii'}\bar{\sigma}_s^{jj'}\bar{\sigma}_s^{pk} W_s^{(n)i'} W_s^{(n)j'}\,ds \xrightarrow{P} \frac{1}{2}\int_0^t c_s^{ij}\sigma_s^{pk}\,ds,$$

(e) $$n[M^{npij}(\bar{\sigma}), W^k]_t = \sum_{i'j'} n \int_0^t \bar{\sigma}_s^{ii'}\bar{\sigma}_s^{jj'}\bar{\sigma}_s^{pk} \int_{n(s)}^s W_r^{(n)i'}\,dW_r^{j'}\,ds \xrightarrow{P} 0.$$

We finish the proof. $\square$

LEMMA 7.7. *Let $X^1, X^2$ be two continuous local martingales with their quadratic variation process $C_t^{ij} = [X^i, X^j]_t = \int_0^t c_s^{ij}\, ds$, and $X_s^{i(n)} = X_s^i - X_{n(s)}^i$. Let $A_t$ and $\bar{A}_t$ be two adapted continuous processes of finite variation, with $A_t = \int_0^t a_s\, ds$ and $\bar{A}_t = \int_0^t \bar{a}_s\, ds$. If $\int_0^1 (a_s)^2\, ds$, $\int_0^1 (\bar{a}_s)^2\, ds$ and $\int_0^1 \|c_s\|^4\, ds$ are bounded by a constant $\alpha$, where $c_t = (c_t^{ij})_{2\times 2}$, then*

$$n \int_0^t X_s^{1(n)} X_s^{2(n)}\, dA_s \xrightarrow{L^2} \tfrac{1}{2} \int_0^t c_s^{12} a_s\, ds, \tag{80}$$

$$n \int_0^t X_s^{1(n)} A_s^{(n)}\, dX_s^2 \xrightarrow{L^2} 0, \tag{81}$$

$$n \int_0^t X_s^{1(n)} A_s^{(n)}\, d\bar{A}_s \xrightarrow{L^2} 0, \tag{82}$$

$$n \int_0^t \int_{n(s)}^s X_r^{1(n)}\, dX_r^2\, dA_s \xrightarrow{L^2} 0, \tag{83}$$

$$n \int_0^t \int_{n(s)}^s A_r^{(n)}\, dX_r^1\, dX_s^2 \xrightarrow{L^2} 0, \tag{84}$$

$$n \int_0^t \int_{n(s)}^s A_r^{(n)}\, dX_r^1\, d\bar{A}_s \xrightarrow{L^2} 0. \tag{85}$$

PROOF. The Burkholder inequality and the Cauchy inequality are used frequently throughout the proof, and the constant $\kappa$ changes from line to line. Let $Y_s^{(n)} = \int_{n(s)}^s, X_r^{1(n)}\, dX_r^2$, $G_t^n(a) = n \int_0^t Y_s^{(n)} a_s\, ds$, $Z^i = \int_{(i-1)/n}^{i/n} Y_s^{(n)}\, ds$. For (83), we need to show that $G_t^n(a) \xrightarrow{L^2} 0$. By (64), we have

$$\begin{aligned} \mathbb{E} \int_0^1 (X_s^{1(n)})^4\, ds &\le \kappa \mathbb{E} \int_0^1 \left( \int_{n(s)}^s c_r^{11}\, dr \right)^2 ds \\ &\le \kappa/n^2 \mathbb{E} \int_0^1 (c_s^{11})^2\, ds \le \kappa\alpha/n^2. \end{aligned} \tag{86}$$

Then by Cauchy's inequality, (64) and (86),

$$\begin{aligned} \mathbb{E} \int_0^1 (Y_s^{(n)})^2\, ds &\le \kappa \mathbb{E} \int_0^1 \int_{n(s)}^s (X_r^{1(n)})^2 c_r^{22}\, dr\, ds \\ &\le \kappa \left( \mathbb{E} \int_0^1 \int_{n(s)}^s (X_r^{1(n)})^4\, dr\, ds \mathbb{E} \int_0^1 \int_{n(s)}^s (c_r^{22})^4\, dr\, ds \right)^{1/2} \\ &\le \kappa/n \left( \mathbb{E} \int_0^1 (X_s^{1(n)})^4\, ds \mathbb{E} \int_0^1 (c_s^{22})^4\, ds \right)^{1/2} \\ &\le \sqrt{\kappa}\alpha/n^2. \end{aligned}$$

Since $Y_s^{(n)}$ is a martingale with mean zero, $\mathbb{E}(Z^i Z^j) = 0$ for $i \neq j$. Again by the Cauchy inequality, $(Z^i)^2 \leq \int_{(i-1)/n}^{i/n} (Y_s^{(n)})^2\, ds/n$. By (87),

$$\mathbb{E}(G_t^n(1))^2 \leq n^2 \sum_{i=1}^{n} \mathbb{E}(Z^i)^2 \leq n \int_0^1 \mathbb{E}(Y_s^{(n)})^2\, ds \leq \kappa\alpha/n.$$

It follows that $G_t^n(1) \xrightarrow{L^2} 0$, which implies that $G_t^n(a) \xrightarrow{L^2} 0$ for a step function $a$. In the general case there is a step function $b$ such that $\mathbb{E}\int_0^1 (a_s - b_s)^2\, ds \to 0$. Since $G_t^n(b_s) \xrightarrow{L^2} 0$ and

$$\begin{aligned}\mathbb{E}\Big(\sup_{0\leq t\leq 1} |G_t^n(a) - G_t^n(b)|\Big) &\leq n\mathbb{E}\left(\int_0^1 |Y_s^{(n)}(a_s - b_s)|\, ds\right)\\ &\leq n\left(\mathbb{E}\int_0^1 (Y_s^{(n)})^2\, ds \mathbb{E}\int_0^1 (a_s - b_s)^2\, ds\right)^{1/2}\\ &\leq \kappa\sqrt{\alpha}\left(\mathbb{E}\int_0^1 (a_s - b_s)^2\, ds\right)^{1/2},\end{aligned}$$

$G_t^n(a) \xrightarrow{L^2} 0$ for any function $a$. So we have (83). Equation (80) follows from (83), (63) and the following integration by parts formula:

$$X_s^{1(n)} X_s^{2(n)} = \int_{n(s)}^{s} X_r^{1(n)}\, dX_r^2 + \int_{n(s)}^{s} X_r^{2(n)}\, dX_r^1 + C_s^{12} - C_{n(s)}^{12}.$$

Now let $\gamma_n = \sup_{1\leq i\leq n} \int_{(i-1)/n}^{i/n} (a_r)^2\, dr$ Since $\gamma_n \leq \alpha$ and $\gamma_n \to 0$, then $\mathbb{E}[\gamma_n] \to 0$ and $\mathbb{E}[\gamma_n]^2 \to 0$. By the Cauchy inequality, $n(A_s^{(n)})^2 \leq \gamma_n$. By (86), we have the following:

(a) $\mathbb{E}\int_0^1 (X_s^{1(n)})^2\, ds \leq \sqrt{\kappa\alpha}/n$,

(b) $\mathbb{E}\left(n\int_0^t X_s^{1(n)} A_s^{(n)}\, d\bar{A}_s\right)^2$
$\leq \mathbb{E}\int_0^1 n(X_s^{1(n)})^2\, ds\, \mathbb{E}\int_0^1 n(A_s^{(n)})^2 (c_s^{22})^2\, ds \leq \alpha\sqrt{\kappa\alpha}\mathbb{E}[\gamma_n]$,

(c) $\mathbb{E}\left(n\int_0^t X_s^{1(n)} A_s^{(n)}\, dX_s^2\right)^2$
$\leq \kappa\mathbb{E}\left(\int_0^t n^2 (X_s^{1(n)})^2 (A_s^{(n)})^2 c_s^{22}\, ds\right)$
$\leq \kappa\left(\mathbb{E}\int_0^1 n^2 (X_s^{1(n)})^4\, ds\right)^{1/2} \left(\mathbb{E}\int_0^1 n^2 (A_s^{(n)})^4 (c_s^{22})^2\, ds\right)^{1/2}$
$\leq \kappa\alpha\sqrt{\mathbb{E}\gamma_n^2}$.

Thus, (81) and (82) follow.

Next, we will prove (84) and (85). Let $H_s^n = \int_{n(s)}^s A_r^{(n)}\, dX_r^1$. Since

$$\mathbb{E}(nH_s^n)^4 \le \kappa n^4 \mathbb{E}\biggl(\int_{n(s)}^s (A_r^{(n)})^2 c_r^{11}\, dr\biggr)^2 \le \kappa n^2 \mathbb{E}\biggl(\gamma_n \int_{n(s)}^s c_r^{11}\, dr\biggr)^2,$$

by (64), we have

$$\begin{aligned}\mathbb{E}\int_0^1 (nH_s^n)^4\, ds &\le \kappa n^2 \mathbb{E}\biggl[\gamma_n^2 \int_0^1 \biggl(\int_{n(s)}^s c_r^{11}\, dr\biggr)^2 ds\biggr]\\ &\le \kappa \mathbb{E}\biggl[\gamma_n^2 \int_0^1 (c_r^{11})^2\, dr\biggr] \le \kappa\alpha \mathbb{E}[\gamma_n^2].\end{aligned}$$

Then

$$\mathbb{E}\biggl(n\int_0^t \int_{n(s)}^s A_r^{(n)} dX_r^1\, dX_s^2\biggr)^2 \le \kappa \mathbb{E}\int_0^1 (nH_s^n)^2 c_s^{22}\, ds \le \kappa\alpha\biggl(\mathbb{E}\int_0^1 (nH_s^n)^4\, ds\biggr)^{1/2},$$

$$\mathbb{E}\biggl(n\int_0^t \int_{n(s)}^s A_r^{(n)} dX_r^1\, d\bar{A}_s\biggr)^2 \le \alpha \mathbb{E}\int_0^1 (nH_s^n)^2\, ds \le \alpha\biggl(\mathbb{E}\int_0^1 (nH_s^n)^4\, ds\biggr)^{1/2},$$

from which we get (84) and (85), since $\mathbb{E}[\gamma_n^2] \to 0$. $\square$

Department of Mathematics
P.O. Box 118105
University of Florida
Gainesville, Florida 32611-8105
USA
e-mail: yan@math.ufl.edu